\documentclass[graybox,footinfo,dvipsnames]{svmult}

\smartqed
\usepackage{mathptmx}       
\usepackage{helvet}         
\usepackage{courier}        
\usepackage{type1cm}        
\usepackage{graphicx}       

\usepackage{array,colortbl}
\usepackage{amsmath,amsfonts,amssymb,bm} 
\DeclareFontFamily{U}{mathx}{\hyphenchar\font45}
\DeclareFontShape{U}{mathx}{m}{n}{
      <5> <6> <7> <8> <9> <10>
      <10.95> <12> <14.4> <17.28> <20.74> <24.88>
      mathx10
      }{}
\DeclareSymbolFont{mathx}{U}{mathx}{m}{n}
\DeclareFontSubstitution{U}{mathx}{m}{n}
\DeclareMathAccent{\widecheck}      {0}{mathx}{"71}

\renewcommand{\email}[1]{\emailname: #1} 

\usepackage{mathptmx}       
\usepackage{helvet}         
\usepackage{courier}        

\usepackage{makeidx}         
\usepackage{graphicx}        
\usepackage[bottom]{footmisc}

\usepackage{latexsym}
\usepackage{amsmath}
\usepackage{amsfonts}
\usepackage{amssymb}
\usepackage{bm}

\usepackage{url}
\usepackage{algorithm}
\usepackage{algorithmic}
\usepackage[misc,geometry]{ifsym}

\spdefaulttheorem{assumption}{Assumption}{\upshape \bfseries}{\itshape}
\spdefaulttheorem{algo}{Algorithm}{\upshape \bfseries}{\itshape}

\newcommand{\R}{{\mathbb{R}}} 

\DeclareSymbolFont{bbold}{U}{bbold}{m}{n}
\DeclareSymbolFontAlphabet{\mathbbold}{bbold}

\providecommand{\argmax}{\operatorname*{argmax}}

\DeclareSymbolFont{bbold}{U}{bbold}{m}{n}
\DeclareSymbolFontAlphabet{\mathbbold}{bbold}

\usepackage{microtype} 

\usepackage[colorlinks=true,linkcolor=black,citecolor=black,urlcolor=black]{hyperref}
\usepackage{bookmark}
\makeatletter 
  \providecommand*{\toclevel@author}{999}
  \providecommand*{\toclevel@title}{0}
\makeatother

\usepackage{mathtools}
\usepackage{pgfplots}
\pgfplotsset{compat=newest}
\usepackage{tikz}
\usetikzlibrary{plotmarks}
\usetikzlibrary{calc}
\usepackage[utf8]{inputenc}
\newcommand{\expec}{\mathbb{E}}
\newlength{\figurewidth}
\newlength{\figureheight}
\newcommand{\logLogSlopeTriangleIncrease}[5]
{

    \pgfplotsextra
    {
        \pgfkeysgetvalue{/pgfplots/xmin}{\xmin}
        \pgfkeysgetvalue{/pgfplots/xmax}{\xmax}
        \pgfkeysgetvalue{/pgfplots/ymin}{\ymin}
        \pgfkeysgetvalue{/pgfplots/ymax}{\ymax}
        \pgfmathsetmacro{\xArel}{#1}
        \pgfmathsetmacro{\yArel}{#3}
        \pgfmathsetmacro{\xBrel}{#1-#2}
        \pgfmathsetmacro{\yBrel}{\yArel}
        \pgfmathsetmacro{\xCrel}{\xArel}
        \pgfmathsetmacro{\lnxB}{\xmin*(1-(#1-#2))+\xmax*(#1-#2)}
        \pgfmathsetmacro{\lnxA}{\xmin*(1-#1)+\xmax*#1}
        \pgfmathsetmacro{\lnyA}{\ymin*(1-#3)+\ymax*#3}
        \pgfmathsetmacro{\lnyC}{\lnyA+#4*(\lnxA-\lnxB)}
        \pgfmathsetmacro{\yCrel}{\lnyC-\ymin)/(\ymax-\ymin)}
        \coordinate (A) at (rel axis cs:\xArel,\yArel);
        \coordinate (B) at (rel axis cs:\xBrel,\yBrel);
        \coordinate (C) at (rel axis cs:\xCrel,\yCrel);

        \draw[#5]   (A)-- node[pos=0.5,anchor=north] {1}
                    (B)-- 
                    (C)-- node[pos=0.5,anchor=west] {#4}
                    cycle;
    }
}

\newcommand{\logLogSlopeTriangle}[5]
{
\pgfplotsextra
    {
        \pgfkeysgetvalue{/pgfplots/xmin}{\xmin}
        \pgfkeysgetvalue{/pgfplots/xmax}{\xmax}
        \pgfkeysgetvalue{/pgfplots/ymin}{\ymin}
        \pgfkeysgetvalue{/pgfplots/ymax}{\ymax}
        \pgfmathsetmacro{\xArel}{#1}
        \pgfmathsetmacro{\yArel}{#3}
        \pgfmathsetmacro{\xBrel}{#1-#2}
        \pgfmathsetmacro{\yBrel}{\yArel}
        \pgfmathsetmacro{\xCrel}{\xArel}
        \pgfmathsetmacro{\lnxB}{\xmin*(1-(#1-#2))+\xmax*(#1-#2)}
        \pgfmathsetmacro{\lnxA}{\xmin*(1-#1)+\xmax*#1}
        \pgfmathsetmacro{\lnyA}{\ymin*(1-#3)+\ymax*#3}
        \pgfmathsetmacro{\lnyC}{\lnyA-#4*(\lnxA-\lnxB)}
        \pgfmathsetmacro{\yCrel}{\lnyC-\ymin)/(\ymax-\ymin)} 
        \coordinate (A) at (rel axis cs:\xArel,\yArel);
        \coordinate (B) at (rel axis cs:\xBrel,\yBrel);
        \coordinate (C) at (rel axis cs:\xCrel,\yCrel);

        \draw[#5]   (A)-- node[pos=0.5,anchor=north] {1}
                    (B)-- 
                    (C)-- node[pos=0.5,anchor=east] {#4}
                    cycle;
    }
}

\newcommand{\add}{}
\newcommand{\change}{}
\newcommand{\remove}[1]{}

\begin{document}

\title*{A Multilevel Monte Carlo Asymptotic-Preserving Particle Method for Kinetic Equations in the Diffusion Limit}
\titlerunning{MLMC AP Particle Method for Kinetic Equations in the Diffusion Limit}
\author{Emil L{\o}vbak \and Giovanni Samaey \and Stefan Vandewalle}
\institute{
Emil L{\o}vbak \and Giovanni Samaey \and Stefan Vandewalle
\at KU Leuven, Department of Computer Science, NUMA Section, Celestijnenlaan 200A box 2402, 3001 Leuven, Belgium
\email{emil.loevbak@cs.kuleuven.be}, \email{giovanni.samaey@cs.kuleuven.be}, \email{stefan.vandewalle@cs.kuleuven.be}
}
\maketitle

\abstract{We propose a multilevel Monte Carlo method for a particle-based asymptotic-preserving scheme for kinetic equations. Kinetic equations model transport and collision\remove{s} of particles in a position-velocity phase-space. With a diffusive scaling, the kinetic equation converges to an advection-diffusion equation in the limit of zero mean free path. Classical particle-based techniques suffer from a strict time-step restriction to maintain stability in this limit. Asymptotic-preserving schemes provide a solution to this time step restriction, but introduce a first-order error in the time step size. We demonstrate how the multilevel Monte Carlo method can be used as a bias reduction technique to perform accurate simulations in the diffusive regime, while leveraging the reduced simulation cost given by the asymptotic-preserving scheme. We describe how to achieve the necessary correlation between simulation paths at different levels and demonstrate the potential of the approach via \change{numerical} experiments.}

\section{Introduction}
\label{mlmcapsec_intro}

Kinetic equations, modeling particle behavior in a position-velocity phase space, occur in many domains. Examples are plasma physics~\cite{Birdsall2004}, bacterial chemotaxis~\cite{Rousset2011c} and computational fluid dynamics~\cite{Pope1981}. 
Many of these applications exhibit a strong time-scale separation, leading to an unacceptably high simulation cost~\cite{Cercignani1988}. However, one typically is only interested in computing the evolution of some macroscopic quantities of interest. These are usually some moments of the particle distribution, which can be computed as averages over velocity space. The time-scale at which these quantities of interest change is often much slower than the time-scale governing the particle dynamics. The nature of the macroscopic dynamics depends on the scaling of the problem, which can be either \emph{hyperbolic} or \emph{diffusive} \cite{Dimarco2014}. 

The model problem\remove{ considered} in this work is a one-dimensional kinetic equation of the form
\begin{equation}
\label{mlmcapeq_kinetic}
\partial_t f(x,v,t) + v\partial_x f(x,v,t) = Q\left(f(x,v,t)\right),
\end{equation}
where $f(x,v,t)$ represents the \change{distribution of particles} as a function of position $x\in\R$ and velocity $v\in\R$ as it evolves in time \change{$t\in\R^+$}. The left-hand side of\remove{ equation}~\eqref{mlmcapeq_kinetic} represents transport, while $Q(f(x,v,t))$ is a collision operator that results in discontinuous velocity changes. As the collision operator, we take the BGK model \cite{Bhatnagar1954}, which represents linear relaxation to an equilibrium distribution that only depends on the \change{particle} density
\begin{equation}
\rho(x,t) = \int f(x,v,t)dv.
\end{equation}

We introduce a parameter $\epsilon$ that represents the mean free path. When decreasing $\epsilon$, the average time between collisions decreases. In this paper, we consider the diffusive scaling. In that case, we simultaneously increase the time scale at which we observe \add{the} evolution of the particle distribution, arriving at
\begin{equation}
\label{mlmcapeq_kineticdiffusive}
\epsilon\partial_t f(x,v,t) + v\partial_x f(x,v,t) = \frac{1}{\epsilon}\left(\add{\mathcal{M}(v)}\rho(x,t)-f(x,v,t)\right),
\end{equation}
\add{with $\mathcal{M}(v)$ the particles' steady state velocity distribution.} It \change{has been} shown that when taking the limit $\epsilon \to 0$, the behavior of equations of the form \eqref{mlmcapeq_kineticdiffusive} is fully described by the diffusion equation \cite{Lapeyre2003}
\begin{equation}
\label{mlmcapeq_heat}
\partial_t \rho(x,t) = \partial_{xx} \rho(x,t).
\end{equation}

\change{Kinetic equations can be simulated} with deterministic methods, solving the partial differential equation (PDE) that describes \add{the} evolution of the particle distribution in \add{the} position-velocity phase space. Alternatively, one can use stochastic methods that simulate a large number of particle trajectories. Deterministic methods 
become prohibitively expensive for higher dimensional applications. Particle-based methods do not suffer from this curse of dimensionality, at the expense of introducing a statistical error in the computed solution. 
The issue of time-scale separation is present in both deterministic and stochastic methods.

One way to avoid the issue of time-scale separation is through the use of asymptotic-preserving methods, which aim at reproducing a scheme for the limiting macroscopic equation in the limit of infinite time-scale separation. For deterministic discretization methods, there is a long line of such methods. We refer to~\cite{Bennoune2008, Boscarino2013, Buet2007, Crouseilles2011, Dimarco2012, Gosse2002, Jin1999, Jin1998, Jin2000, Klar1998, Klar1999, Larsen1974, Lemou2008, Naldi2000} as a representative sample of such methods in the diffusive scaling. The recent review paper \cite{Dimarco2014} contains an overview of the state of the art on asymptotic-preserving methods for kinetic equations, and ample additional references. In the particle-based setting, only a few asymptotic-preserving methods have been developed, mostly in the hyperbolic scaling \cite{Degond2011, Dimarco2008, Dimarco2010, Pareschi1999, Pareschi2001, Pareschi2005}.  In the diffusive scaling, there are only \remove{of} two works~\cite{Crestetto2018,Dimarco2018} so far, to the best of our knowledge. Both methods 
avoid the time step restrictions caused by fast problem time-scales, at the expense of introducing a bias, which is of order one in the time step size. 

The goal of the present paper is to combine the asymptotic-preserving scheme in~\cite{Dimarco2018} with the multilevel Monte Carlo method. Given a fixed computational budget, a trade-off typically has to be made between a small bias and a low variance. The former can be obtained by reducing the time step, the latter by simulating many trajectories with large time steps. The core idea behind the multilevel Monte Carlo method~\cite{Giles2008} is to reduce computational cost, by combining estimates computed with different time step sizes. The multilevel Monte Carlo method\add{,}\remove{ was} originally developed in the context of stochastic processes\add{,}\remove{ and} has been applied to problems across many fields, for example, finance~\cite{Giles2008} and \change{biochemistry}~\cite{Anderson2011}. The method has\remove{ also} successfully been applied to simulating large PDE's with random coefficients~\cite{Cliffe2011}. Recent work has also used multilevel Monte Carlo methods in an optimization context~\cite{VanBarel2019}.

The remainder of this paper is organized as follows. In Section~\ref{mlmcapsec_kinetic}, we describe the model kinetic equation on which we will demonstrate our approach, as well as the asymptotic-preserving Monte Carlo scheme that was introduced in \cite{Dimarco2018}. In Section~\ref{mlmcapsec_mlmc}, we cover the multilevel Monte Carlo method that is the core contribution of this paper. 
In Section~\ref{mlmcapsec_experiments}, we present some preliminary experimental results, demonstrating the properties of the new scheme as well as its computational gain. Finally, in Section~\ref{mlmcapsec_conclusion} we will summarize our main results and mention some possible future extensions.

\section{Model problem and asymptotic-preserving scheme}
\label{mlmcapsec_kinetic}

\subsection{Model equation in the diffusive limit}

The model problem considered in this work is a one-dimensional kinetic equation in the diffusive scaling of the form~\change{\eqref{mlmcapeq_kineticdiffusive}}, which we rewrite as 
\begin{equation}
\label{mlmcapeq_kineticdimless}
\partial_t f(x,v,t) + \dfrac{v}{\epsilon}\partial_x f(x,v,t) = \frac{1}{\epsilon^2}\left(\add{\mathcal{M}(v)}\rho(x,t)-f(x,v,t)\right).
\end{equation}
For ease of exposition, we restrict ourselves to the case of two discrete velocities, $v = \pm 1$. Then, we can write $f_+(x,t)$ and $f_-(x,t)$ to represent the distribution of particles with\change{, respectively,} positive and negative velocities, and $\rho(x,t) = f_+(x,t) + f_-(x,t)$ represents the total density of particles. In this case, equation~\eqref{mlmcapeq_kineticdimless} simplifies to
\begin{equation}
\label{mlmcapeq_GT}
\begin{dcases}
\partial_t f_+(x,t) + \frac{1}{\epsilon} \partial_x f_+(x,t) = \frac{1}{\epsilon^2} \left( \frac{\rho(x,t)}{2} - f_+(x,t) \right) \\
\partial_t f_-(x,t) - \frac{1}{\epsilon} \partial_x f_-(x,t) = \frac{1}{\epsilon^2} \left( \frac{\rho(x,t)}{2} - f_-(x,t) \right)
\end{dcases}.
\end{equation}
Equation~\eqref{mlmcapeq_GT} is also known as the Goldstein-Taylor model, and  can be solved using a particle scheme. \change{For this, we} introduce a time step $\Delta t$ and an ensemble of $P$ particles
\begin{equation}
\left\{\left(X_{p,\Delta t}^n,V_{p,\Delta t}^n\right)\right\}_{p=1}^P.
\end{equation}
The particle state (position and velocity) is represented as $(X,V)$, $p$ is the particle index ($1\le p \le P$), and $n$ represents the time index, i.e., $X_{p,\Delta t}^n\approx X_p(n\Delta t)$. Equation~\eqref{mlmcapeq_GT} is then solved via operator splitting as
\begin{enumerate}
\item \textbf{Transport step.} The position of each particle is updated based on its velocity 
\begin{equation}
\label{mlmceq_transport}
X^{n+1}_{p,\Delta t} = X^n_{p,\Delta t} \change{+} V_{p,\Delta t}^n \Delta t.
\end{equation}
\item \textbf{Collision step.} During collisions, each particle's velocity is updated as:
\begin{equation}
\label{mlmceq_collision}
V_{p,\Delta t}^{n+1} = 
\begin{cases}
\pm 1/\epsilon,&\text{with probability } p_{c,\Delta t} = \Delta t/\epsilon^2 \text{ and equal probability in the sign},\\
V_{p,\Delta t}^n,&\text{otherwise.} 
\end{cases}
\end{equation}
\end{enumerate}
\remove{To maintain stability,} This approximation requires a time step restriction $\Delta t = \mathcal{O}(\epsilon^2)$ \add{as $\epsilon \rightarrow 0$, both to ensure $p_{c,\Delta t} < 1$ in the collision phase, and to keep the increments in the transport phase finite.} \change{This leads} to unacceptably high computational costs \add{for small $\epsilon$}.

\subsection{Asymptotic-preserving Monte Carlo scheme\label{sec:ap_particle}}

Recently, an asymptotic-preserving Monte Carlo scheme was proposed~\cite{Dimarco2018}, based on the simulation of a modified equation
\begin{equation}
\label{mlmcapeq_GTmod}
\begin{dcases}
\partial_t f_+ + \frac{\epsilon}{\epsilon^2+\Delta t} \partial_x f_+ = \frac{\Delta t}{\epsilon^2+\Delta t} \partial_{xx} f_+ + \frac{1}{\epsilon^2+\Delta t} \left( \frac{\rho}{2} - f_+ \right) \\
\partial_t f_- - \frac{\epsilon}{\epsilon^2+\Delta t} \partial_x f_- = \frac{\Delta t}{\epsilon^2+\Delta t} \partial_{xx} f_- + \frac{1}{\epsilon^2+\Delta t} \left( \frac{\rho}{2} - f_- \right)
\end{dcases}.
\end{equation}
In \eqref{mlmcapeq_GTmod} we have dropped the space and time dependency of $f_{\pm}$ and $\rho$, for conciseness. The model given by \eqref{mlmcapeq_GTmod} reduces to \eqref{mlmcapeq_GT} in the limit when $\Delta t$ tends to zero\remove{,} and has an $\mathcal{O}(\Delta t)$ bias. In the limit when $\epsilon$ tends to zero, the equations reduce to \eqref{mlmcapeq_heat}. 

Discretizing this equation, using operator splitting as above, again leads to a Monte Carlo scheme. For each particle $X_p$ and for each time step $n$, one time step now consists of a transport-diffusion and a collision step:
\begin{enumerate}
\item \textbf{Transport-diffusion step.} The position of the particle is updated based on its velocity and a Brownian increment
\begin{equation}
\label{mlmcapeq_transport}
\begin{aligned}
X^{n+1}_{p,\Delta t} &= X^n_{p,\Delta t} \pm \frac{\epsilon}{\epsilon^2+\Delta t} \Delta t + \sqrt{2 \Delta t}\sqrt{\frac{\Delta t}{\epsilon^2 + \Delta t}}\xi^n_p \\
&= X^n_{p,\Delta t} + V^n_{p,\Delta t} \Delta t + \sqrt{2 \Delta t}\sqrt{D_{\Delta t}}\xi^n_p,
\end{aligned}
\end{equation}
in which we have taken $\xi_p^n \sim \mathcal{N}(0,1)$ and introduced a $\Delta t$-dependent velocity $V^n_{p,\Delta t}$ and diffusion coefficient $D_{\Delta t}$:
\begin{equation}\label{mlmcapeq_coefficients}
V^n_{p,\Delta t} = \add{\pm} \frac{\epsilon}{\epsilon^2+\Delta t}, \qquad D_{\Delta t}=\frac{\Delta t}{\epsilon^2 + \Delta t}.
\end{equation}
\item \textbf{Collision step.} During collisions, each particle's velocity is updated as:
\begin{equation}
\label{mlmcapeq_collision}
V_{p,\Delta t}^{n+1} = 
\begin{cases}
\pm \dfrac{\epsilon}{\epsilon^2 + \Delta t},&\text{with probability } \change{p_{c,\Delta t}} = \dfrac{\Delta t}{\epsilon^2+\Delta t} \\
&\text{\phantom{blabla} and equal probability in the sign},\\
V_{p,\Delta t}^n,&\text{otherwise.} 
\end{cases}
\end{equation}
\end{enumerate}
For more details, we refer the reader to~\cite{Dimarco2018}.

\section{Multilevel Monte Carlo method}
\label{mlmcapsec_mlmc}

\subsection{Method and notation}

We want to estimate some quantity of interest $Y$ that is a function of the particle distribution $f(x,v,t)$ at some specific moment $t=t^*$ in time, i.e., we are interested in
\begin{equation}\label{eq:QoI}
Y(t^*) = \expec [F(X(t^*))] = \int\int F(x) f(x,v,t^*)dxdv.
\end{equation}
Note that, in equation~\eqref{eq:QoI}, the function $F$ only depends on the position $x$ and not on velocity. This is a choice we make for notational convenience\remove{,} and is not essential for the method we present.

\change{The classical Monte Carlo estimator $\hat{Y}(t^*)$ for \eqref{eq:QoI} is given by}
\begin{equation}
\label{mlmcapeq_mcestimator}
\hat{Y}(t^*) = \frac{1}{P}\sum_{p=1}^P F(X^N_{p,\Delta t}), \quad t^* = N \Delta t.
\end{equation}
Here, $P$ denotes the number of simulated trajectories, $N$ the number of simulated time steps, $\Delta t$ the time step size\add{,} and $X^N_{p,\Delta t}$ is generated by the time-discretised process \eqref{mlmcapeq_transport}--\eqref{mlmcapeq_collision}. Given a constrained computational budget, a trade-off has to be made when selecting the time step size $\Delta t$. On the one hand, a small time step reduces the bias of the simulation of each sampled trajectory, and thus of the estimated quantity of interest. On the other hand, a large time step reduces the cost per trajectory, which increases the number of trajectories that can be simulated and thus reduces the resulting variance on the estimate. The key idea behind the Multilevel Monte Carlo method~\cite{Giles2008} is to generate a sequence of estimates with varying discretization accuracy and a varying number of realizations. The method achieves the bias of the finest discretization, with the variance of the coarsest discretization.

To apply the multilevel Monte Carlo method, we define a sequence of time step\remove{s} \add{sizes}, denoted by $\Delta t_\ell$ with $\ell = 0\dots L$, with $\ell=L$ denoting the finest level of discretization (smallest time step), and $\ell=0$ the coarsest level. We use a fixed ratio of time steps between subsequent levels, i.e., we set  $\Delta t_{\ell-1}=M\Delta t_\ell$ for some integer $M$. At each level, we simulate a number $P_\ell$ of particle trajectories. An initial coarse estimator with a large number $P_0$ of sample trajectories is  given by
\begin{equation}\label{mlmcapeq_coarse_estimator}
\hat{Y}_0(t^*) = \frac{1}{P_0}\sum_{p=1}^{\change{P_0}} F(X^{N_0}_{p,\Delta t_0}), \quad t^* = \change{N_0} \Delta t_0.
\end{equation}
This initial estimate can be improved upon by a series of difference estimators $\hat{Y}_\ell(t^*)$, $\ell = 1 \dots L$\add{,} of the form
\begin{equation}
\label{mlmcapeq_diffestimator}
\hat{Y}_\ell(t^*) = \frac{1}{P_\ell}\sum_{p=1}^{P_\ell}\left(F(X^{N_\ell}_{p,\Delta t_\ell})-F(X^{N_{\ell-1}}_{p,\Delta t_{\ell-1}})\right),
\end{equation}
with $N_\ell\Delta t_\ell=t^*$, for each value of $\ell$, and $P_\ell$ the number of correlated sample trajectories at each level.
The estimators~\eqref{mlmcapeq_diffestimator} estimate the bias induced by sampling with a simulation time step size $\Delta t_{\ell-1}$ by comparing the sample results with a simulation using a time step size $\Delta t_\ell$. The estimators \eqref{mlmcapeq_coarse_estimator}-\eqref{mlmcapeq_diffestimator} are then combined into a multilevel Monte Carlo estimator via a telescopic sum,
\begin{equation}
\label{mlmcapeq_telescopic}
\hat{Y}(t^*) = \sum_{\ell=0}^{L} \hat{Y}_\ell(t^*).
\end{equation}

It can easily be seen that the expected value of estimator~\eqref{mlmcapeq_telescopic} is the same as that of estimator~\eqref{mlmcapeq_mcestimator} with the finest time step $\Delta t_L$. If the required number of particles $P_\ell$ at each level decreases sufficiently fast with increasing level $\ell$, the multilevel estimator will result in a reduced computational cost for a given accuracy. For more details on the multilevel Monte Carlo method, we refer to~\cite{Giles2015}.

\subsection{Correlating asymptotic-preserving Monte Carlo simulations}
\label{mlmcapsec_correlation}

\subsubsection{Coupled trajectories and notation} 
\label{mlmcapsec_notation}

The differences in \eqref{mlmcapeq_diffestimator} will only have low variance if the simulated paths $X^{n,m}_{\Delta t_\ell,p}$ and $X^n_{\Delta t_{\ell-1},p}$ are correlated. To achieve this correlation, we will couple the different sources of randomness in the simulation at consecutive levels. In each time step using the asymptotic-preserving particle scheme \eqref{mlmcapeq_transport}--\eqref{mlmcapeq_collision}, there are two  sources of stochastic behavior. On the one hand, a new Brownian increment $\xi_p^n$ is generated for each particle in each transport-diffusion step~\eqref{mlmcapeq_transport}.
On the other hand, in each collision step~\eqref{mlmcapeq_collision}, a fraction of particles randomly get a new velocity $V_p^n$. 

Particle trajectories can be coupled by separately correlating the random numbers used for \change{the individual particles} in the transport-diffusion and collision phase of each time step. To show how this is done, we introduce a pair of simulations spanning a time step\remove{ at level $\ell-1$,} with size $\Delta t_{\ell-1}$: (i) a simulation at level $\ell-1$, using a single time step of size $\Delta t_{\ell-1}$; and (ii) a simulation at level $\ell$, using $M$ time steps of size $\Delta t_\ell$:
\begin{equation}
\label{mlmcapeq_transportdiffusionpair}
\begin{dcases}
X^{n+1}_{p,\Delta t_{\ell-1}} \!\!\!\!= X^{n}_{p,\Delta t_{\ell-1}} \!+ \Delta t_{\ell-1}V_{p,\Delta t_{\ell-1}}^n \!+ \sqrt{2 \Delta t_{\ell-1}} \sqrt{D_{\Delta t_{\ell-1}}} \xi^{n}_{p,\ell-1}, \quad\! \xi^{n}_{p,\ell-1} \sim \mathcal{N}(0,1), \\
X^{n+1,0}_{p,\Delta t_{\ell}} = X^{n,0}_{p,\Delta t_{\ell}} + \sum_{m=1}^M \left( \Delta t_\ell V_{p,\Delta t_\ell}^{n,m} + \sqrt{2 \Delta t_{\ell}} \sqrt{D_{\Delta t_\ell}} \xi^{n,m}_{p,\ell}\right), \qquad\quad\; \xi^{n,m}_{p,\ell} \sim \mathcal{N}(0,1),
\end{dcases}
\end{equation}
 with $m \in \{1,\dots, M\}$ and $X^{n,m}_{p,\Delta t_\ell}\approx X_p(n\Delta t_{\ell-1}+m\Delta t_{\ell})\equiv X_p((nM+m)\Delta t_{\ell})$.  
 
 The key point of the algorithm is to compute the velocities $V^n_{p,\Delta t_{\ell-1}}$ and the Brownian increments $\xi^{n}_{p,\ell-1}$ at level $\ell-1$, based on the \remove{based on the} randomly generated values $\xi^{n,m}_{p,\ell}$ and $V^{n,m}_{p,\Delta t_\ell}$ at level $\ell$, instead of generating these independently.  The main difficulty lies in maximizing the correlation between the velocities and Brownian increments at levels $\ell$ and $\ell-1$ while avoiding the introduction of \change{an extra} bias at level $\ell-1$.
Once the coupled simulation~\eqref{mlmcapeq_transportdiffusionpair} at level $\ell-1$ is performed, we can insert the results in~\eqref{mlmcapeq_diffestimator} to obtain a low-variance difference estimator.  In the next two subsections, we explain how we correlate the Brownian increments during the transport phase (Section~\ref{sec:couple_transport}) and the velocities during the collision phase (Section~\ref{sec:couple_collision}). We present the complete algorithm in Section~\ref{sec:algorithm}.

\subsubsection{Coupling the transport-diffusion phase\label{sec:couple_transport}} 

We first correlate the Brownian increments at levels $\ell$ and $\ell-1$. To this end, we first simulate the stochastic process at level $\ell$, using i.d.d.~increments $\xi^{n,m}_{p,\ell}$. Then, at level $\ell-1$, we\remove{ will} compute the Brownian increments, $\xi^{n}_{p,\ell-1}$, from those at level $\ell$, $\left\{\xi^{n,m}_{p,\ell}\right\}$, ensuring that $\xi^{n}_{p,\ell-1}\sim \mathcal{N}(0,1)$. This condition is clearly satisfied if we define $\xi^{n}_{p,\ell-1}$ as
\begin{equation}
\label{mlmcapeq_xicorr}
\xi^{n}_{p,\ell-1} = \sum_{m=1}^M \frac{\xi^{n,m}_{p,\ell}}{\sqrt{M}}.
\end{equation}
Correlating the simulations in this way means that both levels use the same Brownian path, and differences in the diffusion part of the motion only result from differences in the diffusion coefficients $D_\ell$ and $D_{\ell-1}$ at different levels.

\begin{figure}
\centering
%
%
\definecolor{mycolor1}{rgb}{0.12156862745098,0.466666666666667,0.705882352941177}
\definecolor{mycolor2}{rgb}{1,0.498039215686275,0.0549019607843137}
\begin{tikzpicture}

\begin{axis}[%
width=1.8\figurewidth,
height=0.5\figureheight,
at={(0\figurewidth,0\figureheight)},
scale only axis,
xlabel={Time},
ylabel={Position},
xmin=0,
xmax=10,
ymin=-5,
ymax=1,
axis background/.style={fill=white},
legend entries={{$X^{n}_{p,\Delta t_{\ell}}$},{$X^{n}_{p,\Delta t_{\ell-1}}$}},
legend cell align={left},
legend style={at={(0.03,0.08)}, anchor=south west, draw=white!80.0!black}
]
\addlegendimage{mycolor1, mark=square*, mark size=1}
\addlegendimage{mycolor2, mark=*, mark size=1}
\addplot [semithick, color=mycolor1, mark=square*, mark options={}, mark size=1pt, forget plot]
  table[row sep=crcr]{%
0	0\\
0.2	-0.0884229073368045\\
0.4	0.175180492605981\\
0.6	0.25243589000908\\
0.8	-0.181752229749213\\
1	0.218474835734754\\
1.2	0.347943508695561\\
1.4	0.404938268644316\\
1.6	0.622185198819496\\
1.8	0.732403782657772\\
2	0.335438526543269\\
2.2	0.266990953728434\\
2.4	0.205408912753681\\
2.6	-0.01890677949138\\
2.8	0.69033036977426\\
3	0.321090456745502\\
3.2	0.117096120129971\\
3.4	-0.18311135711372\\
3.6	-0.678202746988402\\
3.8	-0.759258044550535\\
4	-0.874816199941612\\
4.2	-0.229680982501482\\
4.4	-0.334679033178855\\
4.6	-0.783390806676384\\
4.8	-0.107313847625855\\
5	0.413272590199273\\
5.2	0.316453577390658\\
5.4	-0.318599113881461\\
5.6	-0.506070662020094\\
5.8	-0.571821175830221\\
6	-0.455420579551636\\
6.2	-0.565536841243248\\
6.4	-0.378573746682483\\
6.6	-0.213336639530805\\
6.8	-0.7406691684776\\
7	-1.14036458813111\\
7.2	-1.45284235408876\\
7.4	-1.66695703338024\\
7.6	-1.80212353510279\\
7.8	-1.79686612564391\\
8	-3.07407943703199\\
8.2	-3.26677372892507\\
8.4	-2.74291148366623\\
8.6	-3.19267228421021\\
8.8	-2.7989779229455\\
9	-2.65126961941534\\
9.2	-2.66349952322763\\
9.4	-2.5865709347805\\
9.6	-3.24645649084197\\
9.8	-3.28210146529928\\
10	-2.60581830361661\\
};
\addplot [semithick, color=mycolor2, mark=*, mark options={}, mark size=1pt, forget plot]
  table[row sep=crcr]{%
0	0\\
1	0.293114750445606\\
2	0.450038008573871\\
3	0.430788052925638\\
4	-1.1736890945325\\
5	0.554463362953792\\
6	-0.611010824537885\\
7	-1.52995964291682\\
8	-4.12431035366269\\
9	-3.55705145741561\\
10	-3.49607211833995\\
};
\end{axis}

\end{tikzpicture}%
\caption{Correlated diffusion steps with $\epsilon=0.5$, $\Delta t_\ell=0.2$ and $\Delta t_{\ell-1}=1$.\label{mlmcapfig_corrpathdiff}}
\end{figure}
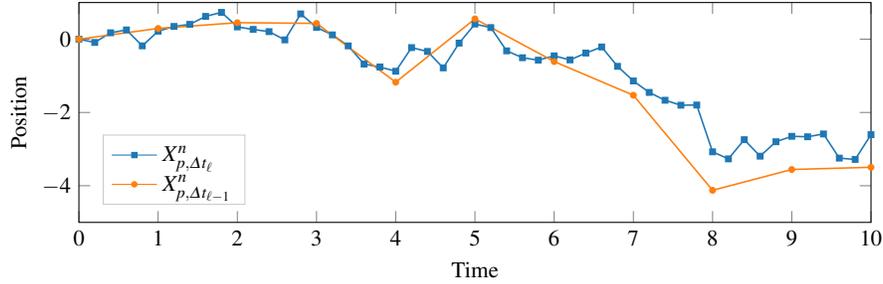

In Figure \ref{mlmcapfig_corrpathdiff}, we show two particle trajectories, containing only diffusion behavior, i.e., \eqref{mlmcapeq_transportdiffusionpair} with $V_{p,\Delta t_{\ell-1}}^n = V_{p,\Delta t_{\ell}}^n = 0$, coupled as described in \eqref{mlmcapeq_xicorr} with $\epsilon=0.5$, $\Delta t_\ell=0.2$ and $M=5$. We observe that the paths have similar behavior, i.e., if the fine simulation tends towards negative values, so does the coarse simulation and vice versa. Still, there is an observable difference between them. This is due to the bias caused by the paths having different diffusion coefficients.

\subsubsection{Coupling the collision phase\label{sec:couple_collision}}

While correlating the Brownian paths is relatively straightforward, the coupling of the velocities in the collision phase is more involved. Since we simulate level $\ell$ first, we have at our disposal the velocities $V^{n,m}_{p,\Delta t_\ell}$ at level $\ell$, which are again i.i.d.  Our goal is to compute the velocities $V^{n}_{p,\Delta t_{\ell-1}}$ at level $\ell-1$ from those at level $\ell$, to maximize correlation, while ensuring that the collision probability and post-collision velocity distribution at level $\ell-1$ are satisfied.  Note that, in the collision phase of the asymptotic-preserving particle scheme~\eqref{mlmcapeq_collision}, both the value of the velocity and the probability of collision depend on the value of the time step $\Delta t$, and therefore depend on the level $\ell$.
 
The computation of $V^{n}_{p,\Delta t_{\ell-1}}$ is done in two steps. First, we will couple the occurrence of a collision at level $\ell-1$ to the occurrence of a collision in one of the $M$ sub-steps of the correlated fine simulation. \change{If} we decide to perform a collision both at level $\ell$ and $\ell-1$, we will correlate the new velocities generated in both simulations.

Let us first consider the simulation at level $\ell$. When simulating the collision step, we decide whether a collision has occurred during a time step of length $\Delta t_\ell$ by drawing a random number $\alpha^{n,m}_{p,\ell} \sim \mathcal{U}([0,1])$ and comparing it to the probability that no collision has occurred in the simulation, $p_{nc,\Delta t_\ell}=1-p_{c,\Delta t_\ell}$, with $p_{c,\Delta t_\ell}$, defined in equation~\eqref{mlmcapeq_collision}. A collision takes place if and only if
\begin{equation}
\label{mlmcapeq_collisioncondition}
\alpha_{p,\ell}^{n,m} \geq p_{nc,\Delta t_\ell} = \frac{\epsilon^2}{\epsilon^2+\Delta t_\ell}.
\end{equation}
Now consider $M$ time steps of length $\Delta t_\ell$. At least one collision has taken place if at least one of the generated $\alpha_{p,\ell}^{n,m}$, $m \in \{1,\dots,M\}$, satisfies \eqref{mlmcapeq_collisioncondition}.\\

\textbf{Deciding upon collision \change{in the coarse simulation}} At level $\ell-1$, we want to use the values $\alpha_{p,\ell}^{n,m}$, $m \in \{1,\dots,M\}$ to compute\remove{d} a uniformly distributed number $\alpha_{p,\ell-1}^n$\add{,} that is correlated with the largest of the generated $\alpha_{p,\ell}^{n,m}$
\begin{equation}
\label{mlmcapeq_maxalpha}
\alpha^{n,\text{max}}_{p,\ell} = \max_m \alpha_{p,\ell}^{n,m}\change{,}
\end{equation}
\add{to compare with the collision probability $p_{nc,\Delta t_{\ell-1}}$.} However, the maximum of a set of uniformly distributed random number is not uniformly distributed. The cumulative density function of $\alpha^{n,\text{max}}_{p,\ell}$ is given by
\begin{equation}\label{eq:max_alpha}
\text{CDF}\left(\alpha^{n,\text{max}}_{p,\ell}\right) = {\left(\alpha^{n,\text{max}}_{p,\ell}\right)}^M.
\end{equation}
Hence, by the inverse transform method, ${\left(\alpha^{n,\text{max}}_{p,\ell}\right)}^M \sim \mathcal{U}([0,1])$.\remove{At the coarse level $\ell-1$, we need a uniformly distributed random number $\alpha_{p,\ell-1}^{n}$, which should be compared to the collision probability $p_{nc,\Delta t_{\ell-1}}$ at level $\ell-1$, equation~(21), using $\ell-1$.} Equation~\eqref{eq:max_alpha} implies that we can define this random number as
\begin{equation}
\label{mlmcapeq_alpha}
\alpha_{p,\ell-1}^n = {\left(\alpha^{n,\text{max}}_{p,\ell}\right)}^M,
\end{equation}
without affecting the simulation statistics at level $\ell-1$.

It is possible to show that, given the relation in \eqref{mlmcapeq_alpha}, a collision can occur in the fine simulation without a collision occurring in the coarse simulation. The inverse, i.e., a collision in the coarse simulation, without a \change{fine simulation collision}, is not possible.\\

\textbf{Choosing a new velocity.} If a collision takes place in both simulations in a given time step $\Delta t_{\ell-1}$, then we set the sign of the velocity of the coarse simulation, at the end of the time step to be equal in sign to the velocity of the last subdividing fine time step for which \eqref{mlmcapeq_collisioncondition} holds,
\begin{equation}
\text{sign}\left(V_{p,\Delta t_{\ell-1}}^{n+1}\right) = \text{sign}\left(V_{p,\Delta t_\ell}^{n,i}\right), \quad i = \argmax_{\quad 1 \leq m \leq M} \left( m \middle| \alpha^{n,m}_{p,\ell} \geq p_{nc,\Delta t_\ell} \right).
\end{equation}
Because the new velocities generated in the fine simulation are i.i.d., we are free to make this selection, without altering the statistics of the coarse simulation. This approach to selecting the sign of $V_{p,\Delta t_{\ell-1}}^{n+1}$ means that the velocities going into the next time step will have the same sign.

\begin{figure}
\centering
%
%
\definecolor{mycolor1}{rgb}{0.12156862745098,0.466666666666667,0.705882352941177}
\definecolor{mycolor2}{rgb}{1,0.498039215686275,0.0549019607843137}
\begin{tikzpicture}

\begin{axis}[%
width=1.8\figurewidth,
height=0.5\figureheight,
at={(0\figurewidth,0\figureheight)},
xlabel={Time},
ylabel={Position},
scale only axis,
xmin=0,
xmax=10,
ymin=0,
ymax=3.5,
axis background/.style={fill=white},
legend entries={{$X^{n}_{p,\Delta t_{\ell}}$},{$X^{n}_{p,\Delta t_{\ell-1}}$}},
legend cell align={left},
legend style={at={(0.03,0.92)}, anchor=north west, draw=white!80.0!black}
]
\addlegendimage{mycolor1, mark=square*, mark size=1}
\addlegendimage{mycolor2, mark=*, mark size=1}
\addplot [semithick, color=mycolor1, forget plot]
  table[row sep=crcr]{%
0	0\\
0.2	0.222222222222222\\
0.4	0.444444444444444\\
0.6	0.666666666666667\\
0.8	0.444444444444445\\
1	0.666666666666667\\
1.2	0.888888888888889\\
1.4	1.11111111111111\\
1.6	0.888888888888889\\
1.8	0.666666666666667\\
2	0.444444444444445\\
2.2	0.666666666666667\\
2.4	0.888888888888889\\
2.6	1.11111111111111\\
2.8	1.33333333333333\\
3	1.55555555555556\\
3.2	1.77777777777778\\
3.4	2\\
3.6	1.77777777777778\\
3.8	1.55555555555556\\
4	1.33333333333333\\
4.2	1.55555555555556\\
4.4	1.77777777777778\\
4.6	2\\
4.8	1.77777777777778\\
5	1.55555555555556\\
5.2	1.33333333333333\\
5.4	1.11111111111111\\
5.6	0.888888888888889\\
5.8	0.666666666666667\\
6	0.888888888888889\\
6.2	1.11111111111111\\
6.4	1.33333333333333\\
6.6	1.11111111111111\\
6.8	1.33333333333333\\
7	1.55555555555556\\
7.2	1.77777777777778\\
7.4	2\\
7.6	2.22222222222222\\
7.8	2.44444444444444\\
8	2.66666666666667\\
8.2	2.88888888888889\\
8.4	3.11111111111111\\
8.6	2.88888888888889\\
8.8	2.66666666666667\\
9	2.44444444444444\\
9.2	2.66666666666667\\
9.4	2.88888888888889\\
9.6	2.66666666666667\\
9.8	2.44444444444444\\
10	2.22222222222222\\
};
\addplot [semithick, color=mycolor2, forget plot]
  table[row sep=crcr]{%
0	0\\
1	0.4\\
2	0.8\\
3	1.2\\
4	1.6\\
5	2\\
6	1.6\\
7	2\\
8	2.4\\
9	2.8\\
10	3.2\\
};
\addplot[only marks, mark=asterisk, mark options={}, mark size=2.5pt, draw=mycolor1] table[row sep=crcr]{%
x	y\\
0.2	0.222222222222222\\
0.4	0.444444444444444\\
0.8	0.444444444444445\\
1.2	0.888888888888889\\
1.4	1.11111111111111\\
1.8	0.666666666666667\\
2	0.444444444444445\\
2.8	1.33333333333333\\
3	1.55555555555556\\
3.2	1.77777777777778\\
3.4	2\\
3.8	1.55555555555556\\
4	1.33333333333333\\
4.6	2\\
5.4	1.11111111111111\\
5.6	0.888888888888889\\
5.8	0.666666666666667\\
6.2	1.11111111111111\\
6.4	1.33333333333333\\
6.6	1.11111111111111\\
7	1.55555555555556\\
7.4	2\\
8	2.66666666666667\\
8.4	3.11111111111111\\
8.8	2.66666666666667\\
9	2.44444444444444\\
9.4	2.88888888888889\\
9.8	2.44444444444444\\
};
\addplot[only marks, mark=asterisk, mark options={}, mark size=2.5pt, draw=mycolor2] table[row sep=crcr]{%
x	y\\
1	0.4\\
2	0.8\\
3	1.2\\
4	1.6\\
5	2\\
6	1.6\\
7	2\\
9	2.8\\
10	3.2\\
};
\addplot[only marks, mark=square*, mark options={}, mark size=1pt, draw=mycolor1, fill=mycolor1] table[row sep=crcr]{%
x	y\\
0	0\\
0.6	0.666666666666667\\
1	0.666666666666667\\
1.6	0.888888888888889\\
2.2	0.666666666666667\\
2.4	0.888888888888889\\
2.6	1.11111111111111\\
3.6	1.77777777777778\\
4.2	1.55555555555556\\
4.4	1.77777777777778\\
4.8	1.77777777777778\\
5	1.55555555555556\\
5.2	1.33333333333333\\
6	0.888888888888889\\
6.8	1.33333333333333\\
7.2	1.77777777777778\\
7.6	2.22222222222222\\
7.8	2.44444444444444\\
8.2	2.88888888888889\\
8.6	2.88888888888889\\
9.2	2.66666666666667\\
9.6	2.66666666666667\\
10	2.22222222222222\\
};
\addplot[only marks, mark=*, mark options={}, mark size=1pt, draw=mycolor2, fill=mycolor2] table[row sep=crcr]{%
0	0\\
8	2.4\\
};
\end{axis}
\end{tikzpicture}%
\caption{Correlated transport steps with $\epsilon=0.5$, $\Delta t_\ell=0.2$ and $\Delta t_{\ell-1}=1$. Stars mark collisions.\label{mlmcapfig_corrpathtransp}}
\end{figure}
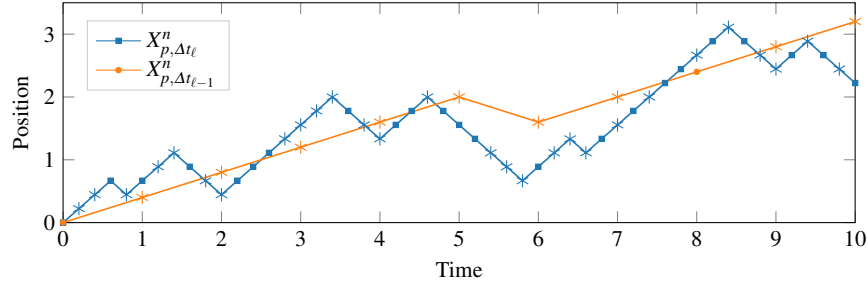

Two particle trajectories without diffusion behavior, i.e., \eqref{mlmcapeq_transportdiffusionpair} with $D_{\Delta t_{\ell-1}} = D_{\Delta t_\ell} = 0$ \change{are} shown in Figure \ref{mlmcapfig_corrpathtransp}. In this figure, a number of interesting phenomena can be observed. First of all, the fact that the particle\add{'s} characteristic velocity is dependent on the time step sizes $\Delta t_{\ell-1}$ and $\Delta t_{\ell}$ results in \change{different slopes in the curves}. This is one source of the bias that we want to estimate using the multilevel Monte Carlo method. \change{Second of all,} the collision probability between the coupled trajectories does not match precisely, as this probability also depends on $\Delta t_{\ell-1}$ and $\Delta t_\ell$. For instance, no collision occurs at $t=8$ in the coarse simulation, while a collision takes place at time $t=7.4$ and $t=8$ in the fine simulation. By coincidence, the new velocity generated at $t=8$ in the fine simulation has the same sign as the \change{coarse simulation velocity}. This mismatch is also part of the bias \remove{that }we wish to estimate.

\subsubsection{The complete algorithm\label{sec:algorithm}}

Combining the correlation of the Brownian increments and velocities results in~Algorithm~\ref{mlmcapalg_correlation}. The correlation of the trajectories \change{can be seen} in Figure \ref{mlmcapfig_corrpath} which shows the particle trajectory given by the sum of the behaviors in Figures \ref{mlmcapfig_corrpathdiff} and \ref{mlmcapfig_corrpathtransp}.
\begin{algorithm}[]
\begin{algorithmic}[1]
\FOR{Each time step $n$}
\FOR{$m = 1 \dots M$}
\STATE Simulate \eqref{mlmcapeq_transport}--\eqref{mlmcapeq_collision} with $\Delta t_\ell$, saving the $\xi^{n,m}_{p,\ell}$, $\alpha^{n,m}_{p,\ell}$ and $V^{n,m}_{p,\Delta t_\ell}$.
\ENDFOR
\STATE Generate $\xi^{n}_{p,\ell-1}$ from the $\xi^{n,m}_{p,\ell}$ according to \eqref{mlmcapeq_xicorr}.
\STATE Generate $\alpha_{p,\ell-1}^n$ from the $\alpha_{p,\ell}^{n,m}$ according to \eqref{mlmcapeq_maxalpha} and \eqref{mlmcapeq_alpha}.
\STATE Set $V_{p,\Delta t_{\ell-1}}^{n+1} = V_{p,\Delta t_{\ell-1}}^{n}$
\IF{$\alpha_{p,\ell-1}^n \geq p_{nc,\Delta t_{\ell-1}}$}
\FOR{$m = 0 \dots M-1$}
\IF{$\alpha_{p,\ell}^{n,m} \geq p_{nc,\Delta t_\ell}$}
\STATE Change the sign of $V_{p,\Delta t_{\ell-1}}^{n+1}$ to be equal to that of $V_{p,\Delta t_{\ell}}^{n,m}$.
\ENDIF
\ENDFOR
\ENDIF
\ENDFOR
\end{algorithmic}
 \caption{Performing correlated simulation steps.\label{mlmcapalg_correlation}}
\end{algorithm}
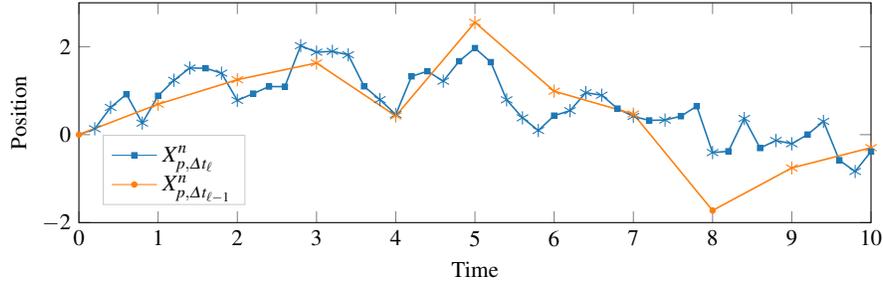
\begin{figure}
\centering
%
%
\definecolor{mycolor1}{rgb}{0.12156862745098,0.466666666666667,0.705882352941177}
\definecolor{mycolor2}{rgb}{1,0.498039215686275,0.0549019607843137}
\begin{tikzpicture}

\begin{axis}[%
width=1.8\figurewidth,
height=0.5\figureheight,
at={(0\figurewidth,0\figureheight)},
xlabel={Time},
ylabel={Position},
scale only axis,
xmin=0,
xmax=10,
ymin=-2,
ymax=3,
axis background/.style={fill=white},
legend entries={{$X^{n}_{p,\Delta t_{\ell}}$},{$X^{n}_{p,\Delta t_{\ell-1}}$}},
legend cell align={left},
legend style={at={(0.03,0.08)}, anchor=south west, draw=white!80.0!black}
]
\addlegendimage{mycolor1, mark=square*, mark size=1}
\addlegendimage{mycolor2, mark=*, mark size=1}
\addplot [semithick, color=mycolor1, forget plot]
  table[row sep=crcr]{%
0	0\\
0.2	0.133799314885418\\
0.4	0.619624937050425\\
0.6	0.919102556675746\\
0.8	0.262692214695232\\
1	0.885141502401421\\
1.2	1.23683239758445\\
1.4	1.51604937975543\\
1.6	1.51107408770838\\
1.8	1.39907044932444\\
2	0.779882970987714\\
2.2	0.9336576203951\\
2.4	1.09429780164257\\
2.6	1.09220433161973\\
2.8	2.02366370310759\\
3	1.87664601230106\\
3.2	1.89487389790775\\
3.4	1.81688864288628\\
3.6	1.09957503078938\\
3.8	0.796297511005021\\
4	0.458517133391721\\
4.2	1.32587457305407\\
4.4	1.44309874459892\\
4.6	1.21660919332362\\
4.8	1.67046393015192\\
5	1.96882814575483\\
5.2	1.64978691072399\\
5.4	0.792511997229651\\
5.6	0.382818226868795\\
5.8	0.0948454908364456\\
6	0.433468309337253\\
6.2	0.545574269867864\\
6.4	0.954759586650851\\
6.6	0.897774471580306\\
6.8	0.592664164855734\\
7	0.415190967424449\\
7.2	0.324935423689017\\
7.4	0.333042966619759\\
7.6	0.420098687119435\\
7.8	0.647578318800531\\
8	-0.407412770365323\\
8.2	-0.377884840036178\\
8.4	0.368199627444877\\
8.6	-0.303783395321323\\
8.8	-0.132311256278836\\
9	-0.206825174970892\\
9.2	0.00316714343903346\\
9.4	0.302317954108386\\
9.6	-0.579789824175304\\
9.8	-0.83765702085484\\
10	-0.383596081394384\\
};
\addplot [semithick, color=mycolor2, forget plot]
  table[row sep=crcr]{%
0	0\\
1	0.693114750445607\\
2	1.25003800857387\\
3	1.63078805292564\\
4	0.426310905467504\\
5	2.55446336295379\\
6	0.988989175462115\\
7	0.470040357083178\\
8	-1.72431035366269\\
9	-0.757051457415612\\
10	-0.296072118339949\\
};
\addplot[only marks, mark=asterisk, mark options={}, mark size=2.5pt, draw=mycolor1] table[row sep=crcr]{%
x	y\\
0.2	0.133799314885418\\
0.4	0.619624937050425\\
0.8	0.262692214695232\\
1.2	1.23683239758445\\
1.4	1.51604937975543\\
1.8	1.39907044932444\\
2	0.779882970987714\\
2.8	2.02366370310759\\
3	1.87664601230106\\
3.2	1.89487389790775\\
3.4	1.81688864288628\\
3.8	0.796297511005021\\
4	0.458517133391721\\
4.6	1.21660919332362\\
5.4	0.792511997229651\\
5.6	0.382818226868795\\
5.8	0.0948454908364456\\
6.2	0.545574269867864\\
6.4	0.954759586650851\\
6.6	0.897774471580306\\
7	0.415190967424449\\
7.4	0.333042966619759\\
8	-0.407412770365323\\
8.4	0.368199627444877\\
8.8	-0.132311256278836\\
9	-0.206825174970892\\
9.4	0.302317954108386\\
9.8	-0.83765702085484\\
};
\addplot[only marks, mark=asterisk, mark options={}, mark size=2.5pt, draw=mycolor2] table[row sep=crcr]{%
x	y\\
1	0.693114750445607\\
2	1.25003800857387\\
3	1.63078805292564\\
4	0.426310905467504\\
5	2.55446336295379\\
6	0.988989175462115\\
7	0.470040357083178\\
9	-0.757051457415612\\
10	-0.296072118339949\\
};
\addplot[only marks, mark=square*, mark options={}, mark size=1pt, draw=mycolor1, fill=mycolor1] table[row sep=crcr]{%
x	y\\
0.6	0.919102556675746\\
1	0.885141502401421\\
1.6	1.51107408770838\\
2.2	0.9336576203951\\
2.4	1.09429780164257\\
2.6	1.09220433161973\\
3.6	1.09957503078938\\
4.2	1.32587457305407\\
4.4	1.44309874459892\\
4.8	1.67046393015192\\
5	1.96882814575483\\
5.2	1.64978691072399\\
6	0.433468309337253\\
6.8	0.592664164855734\\
7.2	0.324935423689017\\
7.6	0.420098687119435\\
7.8	0.647578318800531\\
8.2	-0.377884840036178\\
8.6	-0.303783395321323\\
9.2	0.00316714343903346\\
9.6	-0.579789824175304\\
10	-0.383596081394384\\
};
\addplot[only marks, mark=*, mark options={}, mark size=1pt, draw=mycolor2, fill=mycolor2] table[row sep=crcr]{%
x	y\\
0	0\\
8	-1.72431035366269\\
};
\end{axis}
\end{tikzpicture}%
\caption{Correlated paths steps with $\epsilon=0.5$, $\Delta t_\ell=0.2$ and $\Delta t_{\ell-1}=1$. Stars mark collisions.\label{mlmcapfig_corrpath}}
\end{figure}

\section{Experimental Results}
\label{mlmcapsec_experiments}

We will now demonstrate the viability of the suggested approach through some numerical experiments. We will simulate the model given by \eqref{mlmcapeq_GTmod}, using the multilevel Monte Carlo method to estimate a selected quantity of interest, which is the expected value of the square of the particle position, at \change{$t^*$}. The ensemble of particles is initialized at the origin with equal probability of having a left and right velocity. When discussing results we will replace the full expression for a sample of the quantity of interest, based on an arbitrary particle $p$, $F(X^{N,0}_{\Delta t_\ell,p})$, with the symbol $F_\ell$ to simplify notation.

\subsection{Model correlation behavior}
\label{mlmcapsec_expverification}

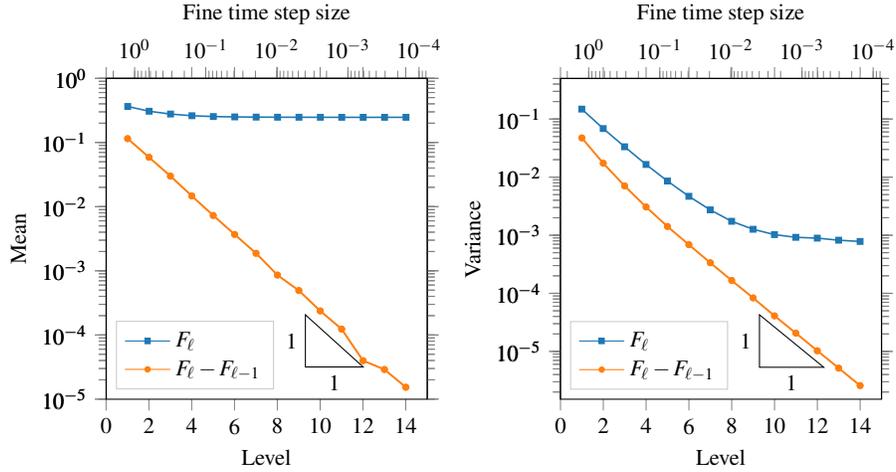
\begin{figure}
\begin{tikzpicture}

\definecolor{color0}{rgb}{0.12156862745098,0.466666666666667,0.705882352941177}
\definecolor{color1}{rgb}{1,0.498039215686275,0.0549019607843137}

\begin{axis}[
title={},
xlabel={Fine time step size},
ylabel={},
xmin=7.62939453125e-05, xmax=2.5,
ymin=10e-06, ymax=1,
xmode=log,
ymode=log,
max space between ticks=20,
xlabel near ticks,
axis x line*=top,
width=\figurewidth,
height=\figureheight,
tick align=outside,
x grid style={white!69.01960784313725!white},
y grid style={white!69.01960784313725!white},
x dir=reverse,
legend entries={{$F_\ell$},{$F_\ell - F_{\ell-1}$}},
legend cell align={left},
legend style={at={(0.03,0.03)}, anchor=south west, draw=white!80.0!black}
]
\addlegendimage{color0, mark=square*, mark size=1}
\addlegendimage{color1, mark=*, mark size=1}
\addplot [semithick, color0, mark=square*, mark size=1, mark options={solid}]
table {%
1.25 0.363268097802
0.625 0.305778636046
0.3125 0.276536394545
0.15625 0.260334031778
0.078125 0.252800667519
0.0390625 0.249472673742
0.01953125 0.247726393937
0.009765625 0.246757812006
0.0048828125 0.246494867005
0.00244140625 0.246195784446
0.001220703125 0.246046872963
0.0006103515625 0.245779917797
0.00030517578125 0.245981342386
0.000152587890625 0.246073464533
};
\addplot [semithick, color1, mark=*, mark size=1, mark options={solid}]
table {%
1.25 0.114791944835
0.625 0.0587506129371
0.3125 0.0299599189201
0.15625 0.0146795430835
0.078125 0.00728449994716
0.0390625 0.00368790241765
0.01953125 0.00187453971133
0.009765625 0.000858235227209
0.0048828125 0.000492204481853
0.00244140625 0.000237160276921
0.001220703125 0.000123324533482
0.0006103515625 3.97741282914e-05
0.00030517578125 2.89664580785e-05
0.000152587890625 1.53182609085e-05
};

\logLogSlopeTriangle{0.62}{-0.18}{0.1}{1}{black}; 

\end{axis}

\begin{axis}[
title={},
xlabel={Level},
ylabel={Mean},
xmin=0, xmax=15,
ymin=10e-06, ymax=1,
ymode=log,max space between ticks=20,
xlabel near ticks,
width=\figurewidth,
height=\figureheight,
tick align=outside,
x grid style={white!69.01960784313725!black},
y grid style={white!69.01960784313725!black},
]
\addplot [semithick, color1, mark=*, mark size=1, mark options={solid}]
table {%
1 0.114791944835
2 0.0587506129371
3 0.0299599189201
4 0.0146795430835
5 0.00728449994716
6 0.00368790241765
7 0.00187453971133
8 0.000858235227209
9 0.000492204481853
10 0.000237160276921
11 0.000123324533482
12 3.97741282914e-05
13 2.89664580785e-05
14 1.53182609085e-05
};

\end{axis}

\end{tikzpicture}
\begin{tikzpicture}

\definecolor{color0}{rgb}{0.12156862745098,0.466666666666667,0.705882352941177}
\definecolor{color1}{rgb}{1,0.498039215686275,0.0549019607843137}

\begin{axis}[
title={},
xlabel=Fine time step size,
ylabel={},
xmin=7.62939453125e-05, xmax=2.5,
ymin=1.5e-06, ymax=0.5,
xmode=log,
ymode=log,max space between ticks=20,
xlabel near ticks,
axis x line*=top,
width=\figurewidth,
height=\figureheight,
tick align=outside,
x grid style={white!69.01960784313725!white},
y grid style={white!69.01960784313725!white},
x dir=reverse,
legend entries={{$F_\ell$},{$F_\ell - F_{\ell-1}$}},
legend cell align={left},
legend style={at={(0.03,0.03)}, anchor=south west, draw=white!80.0!black}
]
\addlegendimage{color0, mark=square*, mark size=1}
\addlegendimage{color1, mark=*, mark size=1}
\addplot [semithick, color0, mark=square*, mark size=1, mark options={solid}]
table {%
1.25 0.148795639866
0.625 0.0687615177621
0.3125 0.0334260928904
0.15625 0.0165664419948
0.078125 0.00859430831669
0.0390625 0.00469210750172
0.01953125 0.00273843967136
0.009765625 0.00174214666174
0.0048828125 0.00126982397042
0.00244140625 0.00102793741574
0.001220703125 0.000922676778289
0.0006103515625 0.000893435765921
0.00030517578125 0.000823147506966
0.000152587890625 0.000783137241272
};
\addplot [semithick, color1, mark=*, mark size=1, mark options={solid}]
table {%
1.25 0.0473493329276
0.625 0.0173490282878
0.3125 0.00706448960306
0.15625 0.003064757657
0.078125 0.00141458967503
0.0390625 0.000687260321554
0.01953125 0.000335842964853
0.009765625 0.000166647955188
0.0048828125 8.34276476792e-05
0.00244140625 4.09447764522e-05
0.001220703125 2.06137230656e-05
0.0006103515625 1.02678456546e-05
0.00030517578125 5.15841712077e-06
0.000152587890625 2.57959641736e-06
};

\logLogSlopeTriangle{0.62}{-0.2}{0.1}{1}{black}; 

\end{axis}

\begin{axis}[
title={},
xlabel={Level},
ylabel={Variance},
xmin=0, xmax=15,
ymin=1.5e-06, ymax=0.5,
ymode=log,max space between ticks=20,
xlabel near ticks,
width=\figurewidth,
height=\figureheight,
tick align=outside,
x grid style={white!69.01960784313725!black},
y grid style={white!69.01960784313725!black},
]
\addplot [semithick, color1, mark=*, mark size=1, mark options={solid}]
table {%
1 0.0473493329276
2 0.0173490282878
3 0.00706448960306
4 0.003064757657
5 0.00141458967503
6 0.000687260321554
7 0.000335842964853
8 0.000166647955188
9 8.34276476792e-05
10 4.09447764522e-05
11 2.06137230656e-05
12 1.02678456546e-05
13 5.15841712077e-06
14 2.57959641736e-06
};

\end{axis}

\end{tikzpicture}
\caption{\label{mlmcapfig_exp1eps10}Mean and variance of the squared particle position for $\epsilon=10$.}
\centering
\end{figure}

In a first test, we \add{set $t^*=5$ and} investigate the variance of the difference estimators~\eqref{mlmcapeq_diffestimator} as a function of the time step $\Delta t_\ell$ (or, equivalently) the level number.  At level $\ell=0$, we set $\Delta t_0 = 2.5$.  All finer levels ($\ell \ge 1$) are defined by setting $\Delta t_\ell=\Delta t_{\ell-1}/M$ with $M=2$. We fix the number of samples per difference estimator at 100~000. For a selection of values of $\epsilon$, we calculate the expected value and variance as a function of $\Delta t_\ell$, for $1\le\ell$. We compute both the variance \change{of} the function samples for a given $\Delta t_\ell$, and the variance of the sampled differences~\eqref{mlmcapeq_diffestimator}, based on the coupled trajectories computed using $\Delta t_{\ell-1}$ and $\Delta t_\ell$. We choose $\epsilon=10$ (Figure~\ref{mlmcapfig_exp1eps10}), $\epsilon=1$ (Figure~\ref{mlmcapfig_exp1eps1}), $\epsilon=0.1$ (Figure~\ref{mlmcapfig_exp1eps01}) and $\epsilon=0.01$ (Figure~\ref{mlmcapfig_exp1eps001}).  \\

\begin{figure}
\begin{tikzpicture}

\definecolor{color0}{rgb}{0.12156862745098,0.466666666666667,0.705882352941177}
\definecolor{color1}{rgb}{1,0.498039215686275,0.0549019607843137}

\begin{axis}[
title={},
xlabel=Fine time step size,
ylabel={},
xmin=1.9073486328125e-05, xmax=2.5,
ymin=4e-06, ymax=20,
xmode=log,
ymode=log,max space between ticks=20,
max space between ticks=20,
xlabel near ticks,
axis x line*=top,
width=\figurewidth,
height=\figureheight,
tick align=outside,
x grid style={white!69.01960784313725!white},
y grid style={white!69.01960784313725!white},
x dir=reverse,
legend entries={{$F_\ell$},{$F_\ell - F_{\ell-1}$}},
legend cell align={left},
legend style={at={(0.03,0.03)}, anchor=south west, draw=white!80.0!black}
]
\addlegendimage{mark=square*, color0, mark size=1}
\addlegendimage{mark=*, color1, mark size=1}
\addplot [semithick, color0, mark=square*, mark size=1, mark options={solid}]
table {%
1.25 7.94001533677
0.625 7.6164871517
0.3125 7.6054460073
0.15625 7.73583949741
0.078125 7.80963650549
0.0390625 7.87716329441
0.01953125 7.94459045039
0.009765625 7.9576559076
0.0048828125 7.98496993476
0.00244140625 8.03304655376
0.001220703125 8.00269556914
0.0006103515625 7.98105317971
0.00030517578125 8.00154231294
0.000152587890625 7.98664204267
7.62939453125e-05 7.9468871424
3.81469726562e-05 8.03647876362
};
\addplot [semithick, color1, mark=*, mark size=1, mark options={solid}]
table {%
1.25 0.575736069864
0.625 0.282648504262
0.3125 0.0269666676531
0.15625 0.115714415131
0.078125 0.131274285832
0.0390625 0.0971562951765
0.01953125 0.0520223085125
0.009765625 0.0247565573833
0.0048828125 0.0107402011465
0.00244140625 0.00435404092613
0.001220703125 0.00418052413651
0.0006103515625 0.000720738464582
0.00030517578125 0.00021838420448
0.000152587890625 0.000497207034614
7.62939453125e-05 0.000245464152265
3.81469726562e-05 7.50383876284e-06
};

\logLogSlopeTriangle{0.62}{-0.25}{0.1}{1}{black}; 

\end{axis}

\begin{axis}[
title={},
xlabel={Level},
ylabel={Mean},
xmin=0, xmax=17,
ymin=4e-06, ymax=20,
ymode=log,max space between ticks=20,
max space between ticks=20,
xlabel near ticks,
width=\figurewidth,
height=\figureheight,
tick align=outside,
x grid style={white!69.01960784313725!black},
y grid style={white!69.01960784313725!black},
]
\addplot [semithick, color1, mark=*, mark size=1, mark options={solid}]
table {%
1 0.575736069864
2 0.282648504262
3 0.0269666676531
4 0.115714415131
5 0.131274285832
6 0.0971562951765
7 0.0520223085125
8 0.0247565573833
9 0.0107402011465
10 0.00435404092613
11 0.00418052413651
12 0.000720738464582
13 0.00021838420448
14 0.000497207034614
15 0.000245464152265
16 7.50383876284e-06
};
\end{axis}

\end{tikzpicture}
\begin{tikzpicture}

\definecolor{color0}{rgb}{0.12156862745098,0.466666666666667,0.705882352941177}
\definecolor{color1}{rgb}{1,0.498039215686275,0.0549019607843137}

\begin{axis}[
title={},
xlabel=Fine time step size,
ylabel={},
xmin=1.9073486328125e-05, xmax=2.5,
ymin=0.00501889080511248, ymax=252.141284796777,
xmode=log,
ymode=log,
max space between ticks=20,
xlabel near ticks,
axis x line*=top,
width=\figurewidth,
height=\figureheight,
tick align=outside,
x grid style={white!69.01960784313725!white},
y grid style={white!69.01960784313725!white},
x dir=reverse,
legend cell align={left},
x dir=reverse,
legend entries={{$F_\ell$},{$F_\ell - F_{\ell-1}$}},
legend style={at={(0.03,0.03)}, anchor=south west, draw=white!80.0!black}
]
\addlegendimage{mark=square*, color0, mark size=1}
\addlegendimage{mark=*, color1, mark size=1}
\addlegendimage{no markers, color2}
\addplot [semithick, color0, mark=square*, mark size=1, mark options={solid}]
table {%
1.25 118.036111864
0.625 100.489903518
0.3125 86.9503554365
0.15625 77.7852146522
0.078125 72.4342294278
0.0390625 68.9299993428
0.01953125 68.4700458316
0.009765625 67.1105737256
0.0048828125 66.8604040883
0.00244140625 66.9648667806
0.001220703125 66.8639968781
0.0006103515625 66.3854315928
0.00030517578125 66.4769142171
0.000152587890625 66.2717249777
7.62939453125e-05 66.0741211735
3.81469726562e-05 66.636515661
};
\addplot [semithick, color1, mark=*, mark size=1, mark options={solid}]
table {%
1.25 22.2781647852
0.625 22.6700174764
0.3125 18.7830829149
0.15625 12.9740475168
0.078125 8.13383212438
0.0390625 4.65351998551
0.01953125 2.49589667182
0.009765625 1.28159831227
0.0048828125 0.65987461745
0.00244140625 0.33611507775
0.001220703125 0.15027075302
0.0006103515625 0.101897965864
0.00030517578125 0.0372707063373
0.000152587890625 0.0323606238426
7.62939453125e-05 0.00845314835635
3.81469726562e-05 0.0108277783303
};

\logLogSlopeTriangle{0.62}{-0.19}{0.1}{1}{black}; 

\end{axis}

\begin{axis}[
title={},
xlabel={Level},
ylabel={Variance},
xmin=0, xmax=17,
ymin=0.00501889080511248, ymax=252.141284796777,
ymode=log,
max space between ticks=20,
xlabel near ticks,
width=\figurewidth,
height=\figureheight,
tick align=outside,
x grid style={white!69.01960784313725!black},
y grid style={white!69.01960784313725!black},
]
\addplot [semithick, color1, mark=*, mark size=1, mark options={solid}]
table {%
1 22.2781647852
2 22.6700174764
3 18.7830829149
4 12.9740475168
5 8.13383212438
6 4.65351998551
7 2.49589667182
8 1.28159831227
9 0.65987461745
10 0.33611507775
11 0.15027075302
12 0.101897965864
13 0.0372707063373
14 0.0323606238426
15 0.00845314835635
16 0.0108277783303
};

\end{axis}

\end{tikzpicture}
\caption{\label{mlmcapfig_exp1eps1} Mean and variance of the squared particle position for $\epsilon=1$.}
\centering
\end{figure}
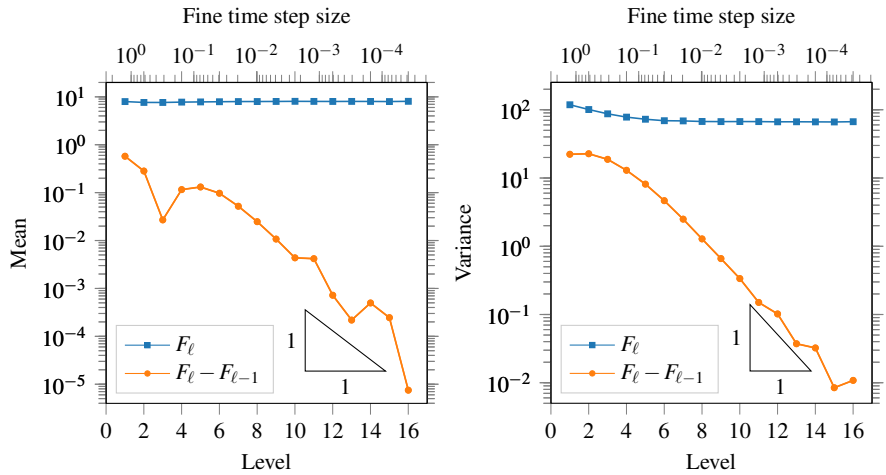
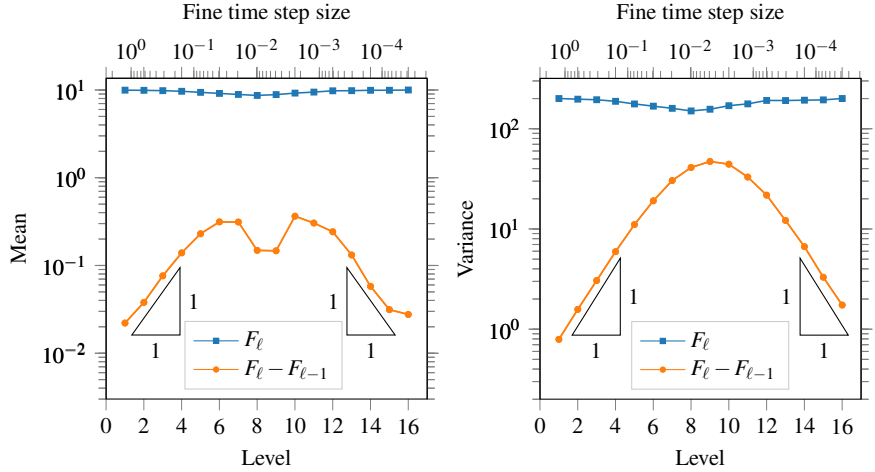
\begin{figure}
\begin{tikzpicture}

\definecolor{color0}{rgb}{0.12156862745098,0.466666666666667,0.705882352941177}
\definecolor{color1}{rgb}{1,0.498039215686275,0.0549019607843137}

\begin{axis}[
title={},
xlabel=Fine time step size,
ylabel={},
xmin=1.9073486328125e-05, xmax=2.5,
ymin=0.003, ymax=13.5970415557483,
xmode=log,
ymode=log,max space between ticks=20,
xlabel near ticks,
axis x line*=top,
width=\figurewidth,
height=\figureheight,
tick align=outside,
x grid style={white!69.01960784313725!white},
y grid style={white!69.01960784313725!white},
x dir=reverse,
legend entries={{$F_\ell$},{$F_\ell - F_{\ell-1}$}},
legend cell align={left},
legend style={at={(0.49,0.03)}, anchor=south, draw=white!80.0!black}
]
\addlegendimage{mark=square*, color0, mark size=1}
\addlegendimage{mark=*, color1, mark size=1}
\addplot [semithick, color0, mark=square*, mark size=1, mark options={solid}]
table {%
1.25 9.96840748524
0.625 9.92701209234
0.3125 9.84183035662
0.15625 9.68277993855
0.078125 9.41395321451
0.0390625 9.1542241289
0.01953125 8.91059774908
0.009765625 8.68166940737
0.0048828125 8.87359245506
0.00244140625 9.23406380224
0.001220703125 9.47984909443
0.0006103515625 9.80526703079
0.00030517578125 9.83856655483
0.000152587890625 9.93754488952
7.62939453125e-05 9.94776834595
3.81469726562e-05 10.0143362108
};
\addplot [semithick, color1, mark=*, mark size=1, mark options={solid}]
table {%
1.25 0.0220890494684
0.625 0.0379730280919
0.3125 0.0762245498094
0.15625 0.139281702302
0.078125 0.229981759774
0.0390625 0.313906904703
0.01953125 0.312839506265
0.009765625 0.14871536068
0.0048828125 0.147006851952
0.00244140625 0.364501735716
0.001220703125 0.304687669943
0.0006103515625 0.242194718957
0.00030517578125 0.131718838343
0.000152587890625 0.057894945633
7.62939453125e-05 0.0313862876206
3.81469726562e-05 0.0276631699499
};

\logLogSlopeTriangleIncrease{0.23}{0.15}{0.2}{1}{black}; 
\logLogSlopeTriangle{0.75}{-0.15}{0.2}{1}{black}; 

\end{axis}

\begin{axis}[
title={},
xlabel={Level},
ylabel={Mean},
xmin=0, xmax=17,
ymin=0.003, ymax=13.5970415557483,
ymode=log,max space between ticks=20,
xlabel near ticks,
width=\figurewidth,
height=\figureheight,
tick align=outside,
x grid style={white!69.01960784313725!black},
y grid style={white!69.01960784313725!black},
]
\addplot [semithick, color1, mark=*, mark size=1, mark options={solid}]
table {%
1 0.0220890494684
2 0.0379730280919
3 0.0762245498094
4 0.139281702302
5 0.229981759774
6 0.313906904703
7 0.312839506265
8 0.14871536068
9 0.147006851952
10 0.364501735716
11 0.304687669943
12 0.242194718957
13 0.131718838343
14 0.057894945633
15 0.0313862876206
16 0.0276631699499
};
\end{axis}

\end{tikzpicture}
\begin{tikzpicture}

\definecolor{color0}{rgb}{0.12156862745098,0.466666666666667,0.705882352941177}
\definecolor{color1}{rgb}{1,0.498039215686275,0.0549019607843137}

\begin{axis}[
title={},
xlabel=Fine time step size,
ylabel={},
xmin=1.9073486328125e-05, xmax=2.5,
ymin=0.2, ymax=317.866980175038,
xmode=log,
ymode=log,
max space between ticks=20,
xlabel near ticks,
axis x line*=top,
width=\figurewidth,
height=\figureheight,
tick align=outside,
x grid style={white!69.01960784313725!white},
y grid style={white!69.01960784313725!white},
x dir=reverse,
legend entries={{$F_\ell$},{$F_\ell - F_{\ell-1}$}},
legend cell align={left},
legend style={at={(0.53,0.03)}, anchor=south, draw=white!80.0!black}
]
\addlegendimage{mark=square*, color0, mark size=1}
\addlegendimage{mark=*, color1, mark size=1}
\addlegendimage{no markers, color2}
\addplot [semithick, color0, mark=square*, mark size=1, mark options={solid}]
table {%
1.25 200.599953961
0.625 197.602810507
0.3125 195.100976482
0.15625 188.13466271
0.078125 177.108466721
0.0390625 168.291019986
0.01953125 160.106371909
0.009765625 150.811056995
0.0048828125 156.82559135
0.00244140625 170.106705416
0.001220703125 177.42238752
0.0006103515625 192.008903765
0.00030517578125 191.599309484
0.000152587890625 193.104990471
7.62939453125e-05 194.436980328
3.81469726562e-05 200.665272706
};
\addplot [semithick, color1, mark=*, mark size=1, mark options={solid}]
table {%
1.25 0.79182562221
0.625 1.57527331776
0.3125 3.05934537554
0.15625 5.93665230953
0.078125 11.0592531206
0.0390625 19.1070602726
0.01953125 30.4314797644
0.009765625 41.055244543
0.0048828125 47.1504126074
0.00244140625 44.3055273615
0.001220703125 32.9694210477
0.0006103515625 21.7525565864
0.00030517578125 12.1576923066
0.000152587890625 6.66458096373
7.62939453125e-05 3.29932453829
3.81469726562e-05 1.74211974356
};

\logLogSlopeTriangleIncrease{0.25}{0.15}{0.2}{1}{black}; 
\logLogSlopeTriangle{0.81}{-0.15}{0.2}{1}{black}; 

\end{axis}

\begin{axis}[
title={},
xlabel={Level},
ylabel={Variance},
xmin=0, xmax=17,
ymin=0.2, ymax=317.866980175038,
ymode=log,max space between ticks=20,
xlabel near ticks,
width=\figurewidth,
height=\figureheight,
tick align=outside,
x grid style={white!69.01960784313725!black},
y grid style={white!69.01960784313725!black},
]
\addplot [semithick, color1, mark=*, mark size=1, mark options={solid}]
table {%
1 0.79182562221
2 1.57527331776
3 3.05934537554
4 5.93665230953
5 11.0592531206
6 19.1070602726
7 30.4314797644
8 41.055244543
9 47.1504126074
10 44.3055273615
11 32.9694210477
12 21.7525565864
13 12.1576923066
14 6.66458096373
15 3.29932453829
16 1.74211974356
};

\end{axis}

\end{tikzpicture}
\caption{\label{mlmcapfig_exp1eps01} Mean and variance of the squared particle position for $\epsilon=0.1$.}
\centering
\end{figure}
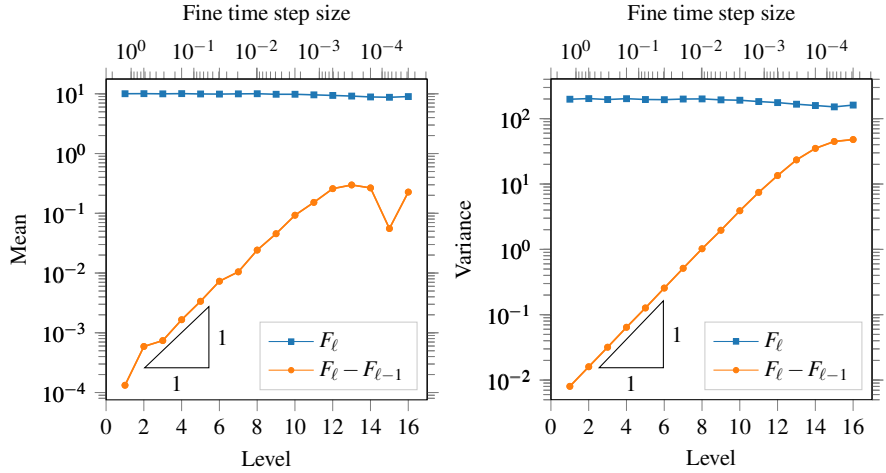
\begin{figure}
\begin{tikzpicture}

\definecolor{color0}{rgb}{0.12156862745098,0.466666666666667,0.705882352941177}
\definecolor{color1}{rgb}{1,0.498039215686275,0.0549019607843137}

\begin{axis}[
title={},
xlabel=Fine time step size,
ylabel={},
xmin=1.9073486328125e-05, xmax=2.5,
ymin=7.54152663311099e-05, ymax=17.6700281591673,
xmode=log,
ymode=log,max space between ticks=20,
xlabel near ticks,
axis x line*=top,
width=\figurewidth,
height=\figureheight,
tick align=outside,
x grid style={white!69.01960784313725!white},
y grid style={white!69.01960784313725!white},
x dir=reverse,
legend entries={{$F_\ell$},{$F_\ell - F_{\ell-1}$}},
legend style={at={(0.97,0.03)}, anchor=south east, draw=white!80.0!black},
legend cell align={left}
]
\addlegendimage{mark=square*, color0, mark size=1}
\addlegendimage{mark=*, color1, mark size=1}
\addplot [semithick, color0, mark=square*, mark size=1, mark options={solid}]
table {%
1.25 10.0332914554
0.625 10.060206185
0.3125 9.98993990723
0.15625 10.0729441195
0.078125 9.9422042789
0.0390625 9.89146703197
0.01953125 9.98076221171
0.009765625 10.0333511082
0.0048828125 9.82068682864
0.00244140625 9.82179382673
0.001220703125 9.60603246538
0.0006103515625 9.39013738906
0.00030517578125 9.16276008958
0.000152587890625 8.86152059235
7.62939453125e-05 8.74415627272
3.81469726562e-05 8.99037781221
};
\addplot [semithick, color1, mark=*, mark size=1, mark options={solid}]
table {%
1.25 0.000132293981173
0.625 0.000590123376795
0.3125 0.000739905339119
0.15625 0.00165452038106
0.078125 0.00336643013131
0.0390625 0.00729403394273
0.01953125 0.0104704337608
0.009765625 0.0242198789932
0.0048828125 0.0455253154734
0.00244140625 0.0927604805724
0.001220703125 0.152236971983
0.0006103515625 0.258347036428
0.00030517578125 0.297942469482
0.000152587890625 0.265981327232
7.62939453125e-05 0.0555163237107
3.81469726562e-05 0.226589289455
};

\logLogSlopeTriangleIncrease{0.32}{0.2}{0.1}{1}{black}; 

\end{axis}

\begin{axis}[
title={},
xlabel={Level},
ylabel={Mean},
xmin=0, xmax=17,
ymin=7.54152663311099e-05, ymax=17.6700281591673,
ymode=log,
max space between ticks=20,
xlabel near ticks,
width=\figurewidth,
height=\figureheight,
tick align=outside,
x grid style={white!69.01960784313725!black},
y grid style={white!69.01960784313725!black},
]
\addplot [semithick, color1, mark=*, mark size=1, mark options={solid}]
table {%
1 0.000132293981173
2 0.000590123376795
3 0.000739905339119
4 0.00165452038106
5 0.00336643013131
6 0.00729403394273
7 0.0104704337608
8 0.0242198789932
9 0.0455253154734
10 0.0927604805724
11 0.152236971983
12 0.258347036428
13 0.297942469482
14 0.265981327232
15 0.0555163237107
16 0.226589289455
};

\end{axis}

\end{tikzpicture}
\begin{tikzpicture}

\definecolor{color0}{rgb}{0.12156862745098,0.466666666666667,0.705882352941177}
\definecolor{color1}{rgb}{1,0.498039215686275,0.0549019607843137}

\begin{axis}[
title={},
xlabel=Fine time step size,
ylabel={},
xmin=1.9073486328125e-05, xmax=2.5,
ymin=0.005, ymax=405.49470514225,
xmode=log,
ymode=log,max space between ticks=20,
xlabel near ticks,
axis x line*=top,
width=\figurewidth,
height=\figureheight,
tick align=outside,
x grid style={white!69.01960784313725!white},
y grid style={white!69.01960784313725!white},
x dir=reverse,
legend style={at={(0.97,0.03)}, anchor=south east, draw=white!80.0!black},
legend cell align={left},
legend entries={{$F_\ell$},{$F_\ell - F_{\ell-1}$}},
]
\addlegendimage{mark=square*, color0, mark size=1}
\addlegendimage{mark=*, color1, mark size=1}
\addlegendimage{no markers, color2}
\addplot [semithick, color0, mark=square*, mark size=1, mark options={solid}]
table {%
1.25 198.862274153
0.625 203.264847785
0.3125 196.820091658
0.15625 202.618168867
0.078125 196.896303242
0.0390625 195.790959567
0.01953125 200.101688317
0.009765625 201.094123104
0.0048828125 194.529087682
0.00244140625 191.928195123
0.001220703125 183.244126641
0.0006103515625 176.509720185
0.00030517578125 166.784253643
0.000152587890625 159.355932779
7.62939453125e-05 151.888375566
3.81469726562e-05 161.727824715
};
\addplot [semithick, color1, mark=*, mark size=1, mark options={solid}]
table {%
1.25 0.00801706694158
0.625 0.0159920064974
0.3125 0.0318247800558
0.15625 0.0642507795671
0.078125 0.126873502589
0.0390625 0.255852942328
0.01953125 0.513408155702
0.009765625 1.02186143721
0.0048828125 1.96156312125
0.00244140625 3.89792605992
0.001220703125 7.45081206185
0.0006103515625 13.4873310251
0.00030517578125 23.4039117139
0.000152587890625 35.1162175386
7.62939453125e-05 44.7572179481
3.81469726562e-05 47.9011176518
};

\logLogSlopeTriangleIncrease{0.35}{0.2}{0.1}{1}{black}; 

\end{axis}

\begin{axis}[
title={},
xlabel={Level},
ylabel={Variance},
xmin=0, xmax=17,
ymin=0.005, ymax=405.49470514225,
ymode=log,
max space between ticks=20,
xlabel near ticks,
width=\figurewidth,
height=\figureheight,
tick align=outside,
x grid style={white!69.01960784313725!black},
y grid style={white!69.01960784313725!black},
]
\addplot [semithick, color1, mark=*, mark size=1, mark options={solid}]
table {%
1 0.00801706694158
2 0.0159920064974
3 0.0318247800558
4 0.0642507795671
5 0.126873502589
6 0.255852942328
7 0.513408155702
8 1.02186143721
9 1.96156312125
10 3.89792605992
11 7.45081206185
12 13.4873310251
13 23.4039117139
14 35.1162175386
15 44.7572179481
16 47.9011176518
};

\end{axis}

\end{tikzpicture}
\caption{\label{mlmcapfig_exp1eps001} Mean and variance of the squared particle position for $\epsilon=0.01$.}
\centering
\end{figure}

\textbf{The regime $\Delta t \ll \epsilon^2$.} In Figures \ref{mlmcapfig_exp1eps10} through \ref{mlmcapfig_exp1eps01}, we see that the slopes of both the mean and variance curves for the differences approach an asymptotic limit $\mathcal{O}\left( \change{\Delta t} \right)$ for $\Delta t \ll \epsilon^2$. This matches the weak convergence order of the Euler-Maruyama scheme, used to simulate the model \eqref{mlmcapeq_transport}--\eqref{mlmcapeq_collision}, as well as the expected behavior from the time step dependent bias in the asymptotic-preserving model. Given this asymptotic geometric convergence, it is possible to apply the complexity theorem in~\cite{Giles2008} to analyze the \add{method's} computational cost and error bounds\remove{ on the method}. This means that \remove{the} existing \remove{general} theory for multilevel Monte Carlo methods~\cite{Giles2015} concerning, e.g. samples per level, convergence criteria and conditions for adding levels, can be applied in this regime.

\textbf{The regime $\Delta t \gg \epsilon^2$.} For time steps $\Delta t \gg \epsilon^2$, however, we see in Figures~\ref{mlmcapfig_exp1eps01} and \ref{mlmcapfig_exp1eps001} that both the mean and the variance curves increase geometrically in terms of increasing level. To explain this perhaps counterintuitive result, we will look at the limit of the modified Goldstein-Taylor model when $\Delta t$ tends to infinity. In this limit, the model \eqref{mlmcapeq_GTmod} converges \add{to} the heat equation:
\begin{equation}
\label{mlmcapeq_GTheat}
\begin{dcases}
\partial_t f_+(x,t) = \partial_{xx} f_+(x,t)\\
\partial_t f_-(x,t) = \partial_{xx} f_-(x,t)
\end{dcases} \Rightarrow \partial_t \rho(x,t) = \partial_{xx} \rho(x,t).
\end{equation}
This means that taking increasingly larger time steps in \eqref{mlmcapeq_GTmod} is equivalent to taking the limit $\epsilon \to 0$. This observation is precisely the asymptotic-preserving property of the particle scheme of Section~\ref{sec:ap_particle}.

\change{The fact that the two limits approach different models can be seen most clearly in Figures~\ref{mlmcapfig_exp1eps10} and~\ref{mlmcapfig_exp1eps01}. In the right hand panel of Figure~\ref{mlmcapfig_exp1eps10} we see that the variance of the individual simulations at level $\ell$ (blue line with squares) changes drastically as a function of $\Delta t_\ell$ in the region where it is of the same order of magnitude as $\epsilon^2$. This is caused by the approximated models for large and small $\Delta t$ having differences in behavior, which are significant enough to be observed when plotted. The scheme thus converges to different equations for the two limits in $\Delta t$. For small $\Delta t$, there is convergence to \eqref{mlmcapeq_GT}. For large $\Delta t$ there is convergence to \eqref{mlmcapeq_heat}. In practice, the size of $\Delta t$ is limited by the simulation time horizon, so it is not possible to get arbitrarily close to \eqref{mlmcapeq_heat} by increasing the time step size, however. This phenomenon also has an effect on the curves in Figure~\ref{mlmcapfig_exp1eps01}. The curves for the mean and variance of the differences $F_\ell - F_{\ell-1}$ (orange lines with dots) decrease for both small and large $\Delta t$, as the model converges to the two limits.}

Combining the observations from the two limits in the time step size gives \remove{us} an intuitive interpretation to the multilevel Monte Carlo method in this setting\change{:} The method can be interpreted as \change{correcting} the result of a \remove{cheap} pure diffusion simulation by decreasing \change{$\Delta t$ to get a good approximation of the transport-diffusion equation that describes the behavior for a given value of $\epsilon$.} The peak of the variance of the differences lies \change{near} $\Delta t \approx \epsilon^2$. This makes sense, as this is the region \change{where} the model parameters $D_{\Delta t}$ and $V^n_{p,\Delta t}$ vary the most in function of $\Delta t$. We also see a dip in the mean of the difference curves in the region of $\Delta t \approx \epsilon^2$. A full analysis of the behavior that occurs in the transition between the asymptotic regimes \change{is} left for \change{future work}.

\subsection{Comparison with classical Monte Carlo}
\label{mlmcapsec_expbenchmark}

The analysis in Section~\ref{mlmcapsec_expverification} demonstrated a fast decay of the variance of the differences for increasingly fine levels in the region where $\Delta t \ll \epsilon^2$. As such, one of the necessary requirements for convergence of the multilevel Monte Carlo method is present in this region. \change{This is, however, not the case in the regime where $\Delta t \gg \epsilon^2$. Here, the variance of the differences increases as the time step is refined. It is therefore highly non-trivial to perform an adequate selection of coarse levels in the regime $\Delta t \gg \epsilon^2$. For the fine levels, a standard multilevel Monte Carlo approach can be applied. We therefore propose two simulation strategies:
\begin{enumerate}
\item A geometric sequence of levels $\Delta t_\ell = \epsilon^2M^{-\ell}$ for $\ell>0$ starting with a coarse simulation time step $\Delta t_0=\epsilon^2$;
\item The same geometric sequence, preceded by a coarse simulation time step $t^*$, i.e., $\Delta t_0=t^*$, $\Delta t_1=\epsilon^2$ and $\Delta t_\ell = \epsilon^2M^{1-\ell}$ for $\ell>1$.
\end{enumerate}
We compare these approaches in the following two sub-sections.}

\subsubsection{Standard MLMC refinement}
We will now compute the quantity of interest described at the beginning of this section to a range of prescribed error tolerances, to verify the reduced computational cost of the multilevel Monte Carlo method. We choose to set $M=2$ and $\epsilon=0.1$\add{, and reduce the time horizon to $t^*=0.5$}. This \remove{value for $\epsilon$} gives us an expensive, but computationally feasible problem. The number of samples per level is derived using the formula~\cite{Giles2015}
\begin{equation}
\left\lceil 2 E^{-2} \sqrt{\frac{V_\ell}{C_\ell}}\left( \sum_{\ell=0}^L \sqrt{V_\ell C_\ell} \right) \right\rceil,
\end{equation}
where $E$ is the desired \add{root} mean square error, $C_\ell$ is the computational cost of the estimator at level $\ell$\add{,} and $V_\ell$ is the estimated variance of the estimator at level $\ell$, i.e., $V_\ell = \mathbb{V}\left[F_\ell - F_{\ell-1} \right]$, where we set $F_{-1} \equiv 0$. The criterion for adding levels and determining convergence are as described in~\cite{Giles2015}. The cost of a sample will be determined relative to the cost of a simulated trajectory with $\Delta t = \epsilon^2$. The results of the simulations for $E$ \change{values 0.1, 0.01 and 0.001} can be found in Tables \ref{mlmcaptab_simulationresults1} through \ref{mlmcaptab_simulationresults3}.

\begin{table}
\caption{\change{Results of the simulation in Section \ref{mlmcapsec_expbenchmark} with a geometric level sequence for} $E=0.1$.\label{mlmcaptab_simulationresults1}}
\centering
\change{
\begin{tabular}{c | c | c | c | r | c | c | c || c}
Level & $\Delta t_\ell$ & $P_\ell$ & $\mathbb{V}\left[ F_\ell \right]$ & \multicolumn{1}{c|}{$\mathbb{E}\left[ F_\ell - F_{\ell-1} \right]$} & \change{$V_\ell$} & $\mathbb{V}[\hat{Y}_\ell]$ & $C_\ell$ & $P_\ell C_\ell$ \\
\hline
0 & $1.00 \times 10^{-2}$ & 1 393 & 1.32 & $8.18 \times 10^{-1}$ & $1.32 \times 10^{0\phantom{-}}$ & $9.45 \times 10^{-4}$ & 1 & 1 393\\
1 & $5.00 \times 10^{-3}$ & 395 & 1.52 & $7.91 \times 10^{-3}$ & $3.58 \times 10^{-1}$ & $9.07 \times 10^{-4}$ & 3 & 1 185\\
2 & $2.50 \times 10^{-3}$ & 296 & 1.59 & $2.18 \times 10^{-2}$  & $4.82 \times 10^{-1}$ & $1.59 \times 10^{-3}$ & 6 & 1 776\\
3 & $1.25 \times 10^{-3}$ & 229 & 2.22 & $-1.48 \times 10^{-2}$  & $3.22 \times 10^{-1}$ & $1.41 \times 10^{-3}$ & 12 & 2 748\\
4 & $6.25 \times 10^{-4}$ & 40  & 1.70 & $1.57 \times 10^{-3}$ & $4.56 \times 10^{-2}$ & $1.14 \times 10^{-3}$ & 24 & 960\\
\hline
$\sum$ & \multicolumn{4}{c}{} & & $6.00 \times 10^{-3}$ & & 8 062
\end{tabular}
}
\end{table}

\begin{table}
\caption{\change{Results of the simulation in Section \ref{mlmcapsec_expbenchmark} with a geometric level sequence for} $E=0.01$.\label{mlmcaptab_simulationresults2}}
\centering
\change{
\begin{tabular}{c | c | c | c | c | c | c | c || c}
Level & $\Delta t_\ell$ & $P_\ell$ & $\mathbb{V}\left[ F_\ell \right]$ & $\mathbb{E}\left[ F_\ell - F_{\ell-1} \right]$ & \change{$V_\ell$} & $\mathbb{V}[\hat{Y}_\ell]$ & $C_\ell$ & $P_\ell C_\ell$ \\
\hline
0 & $1.00 \times 10^{-2}$ & 527 920 & 1.47 & $8.65 \times 10^{-1}$ & $1.47 \times 10^{0\phantom{-}}$ & $2.79 \times 10^{-6}$ & 1   & 527 920\\
1 & $5.00 \times 10^{-3}$ & 165 386  & 1.49 & $1.06 \times 10^{-2}$ & $4.35 \times 10^{-1}$ & $2.63 \times 10^{-6}$ & 3   & 496 158\\
2 & $2.50 \times 10^{-3}$ & 112 208  & 1.59 & $2.98 \times 10^{-2}$ & $3.99 \times 10^{-1}$ & $3.55 \times 10^{-6}$ & 6   & 673 248\\
3 & $1.25 \times 10^{-3}$ & 69 135  & 1.64 & $2.84 \times 10^{-2}$ & $3.01 \times 10^{-1}$ & $4.36 \times 10^{-6}$ & 12  & 829 620\\
4 & $6.25 \times 10^{-4}$ & 39 146  & 1.73 & $2.00 \times 10^{-2}$ & $1.95 \times 10^{-1}$ & $4.98 \times 10^{-6}$ & 24  & 939 504\\
5 & $3.13 \times 10^{-4}$ & 20 670   & 1.76 & $7.53 \times 10^{-3}$  & $1.09 \times 10^{-1}$ & $5.28 \times 10^{-6}$ & 48  & 992 160\\
6 & $1.56 \times 10^{-4}$ & 10 842   & 1.75 & $9.55 \times 10^{-3}$  & $6.14 \times 10^{-2}$ & $5.67 \times 10^{-6}$ & 96  & 1 040 832\\
7 & $7.81 \times 10^{-5}$ & 4 894   & 1.91 & $6.77 \times 10^{-3}$ & $2.42 \times 10^{-2}$ & $4.94 \times 10^{-6}$ & 192 & 939 648\\
8 & $3.91 \times 10^{-5}$ & 3 937   & 1.77 & $2.88 \times 10^{-3}$ & $1.21 \times 10^{-2}$ & $3.08 \times 10^{-6}$  & 384 & 1 511 808\\
9 & $1.95 \times 10^{-5}$ & 2 721   & 1.81 & $2.35 \times 10^{-3}$ & $1.21 \times 10^{-2}$ & $4.46 \times 10^{-6}$ & 768 & 2 089 728\\
10 & $9.75 \times 10^{-6}$ & 40   & 1.47 & $2.00 \times 10^{-3}$ & $4.35 \times 10^{-4}$ & $1.09 \times 10^{-5}$ & 1 536 & 61 440\\
\hline
$\sum$ & \multicolumn{4}{c}{} & & $5.26 \times 10^{-5}$ & & 10 102 066
\end{tabular}
}
\end{table}

\begin{table}
\caption{\change{Results of the simulation in Section \ref{mlmcapsec_expbenchmark} with a geometric level sequence for} $E=0.001$.\label{mlmcaptab_simulationresults3}}
\centering
\noindent\makebox[\textwidth]{
\change{
\begin{tabular}{c | c | c | c | r | c | c | c || c}
Level & $\Delta t_\ell$ & $P_\ell$ & $\mathbb{V}\left[ F_\ell \right]$ & \multicolumn{1}{c|}{$\mathbb{E}\left[ F_\ell - F_{\ell-1} \right]$} & \change{$V_\ell$} & $\mathbb{V}[\hat{Y}_\ell]$ & $C_\ell$ & $P_\ell C_\ell$ \\
\hline
0 & $1.00 \times 10^{-2}$ & 71 593 376 & 1.47 & $8.65 \times 10^{-1}$ & $1.47 \times 10^{0\phantom{-}}$ & $2.06 \times 10^{-8}$ & 1  & 71 593 376\\
1 & $5.00 \times 10^{-3}$ & 22 501 565 & 1.49 & $1.08 \times 10^{-2}$ & $4.36 \times 10^{-1}$ & $1.94 \times 10^{-8}$ & 3   & 67 504 695\\
2 & $2.50 \times 10^{-3}$ & 15 284 042 & 1.57 & $2.84 \times 10^{-2}$ & $4.03 \times 10^{-1}$ & $2.63 \times 10^{-8}$ & 6   & 91 704 252\\
3 & $1.25 \times 10^{-3}$ & 9 372 999  & 1.66 & $2.88 \times 10^{-2}$ & $3.03 \times 10^{-1}$ & $3.23 \times 10^{-8}$ & 12  & 112 475 988\\
4 & $6.25 \times 10^{-4}$ & 5 322 687  & 1.73 & $2.07 \times 10^{-2}$ & $1.95 \times 10^{-1}$ & $3.67 \times 10^{-8}$ & 24  & 127 744 488\\
5 & $3.13 \times 10^{-4}$ & 2 850 794  & 1.77 & $1.24 \times 10^{-2}$ & $1.12 \times 10^{-1}$ & $3.93 \times 10^{-8}$ & 48  & 136 838 112\\
6 & $1.56 \times 10^{-4}$ & 1 480 624  & 1.77 & $6.73 \times 10^{-3}$ & $6.03 \times 10^{-2}$ & $4.07 \times 10^{-8}$ & 96  & 142 139 904\\
7 & $7.81 \times 10^{-5}$ & 749 144    & 1.80 & $3.42 \times 10^{-3}$ & $3.09 \times 10^{-2}$ & $4.13 \times 10^{-8}$ & 192 & 143 835 648\\
8 & $3.91 \times 10^{-5}$ & 382 855    & 1.79 & $1.93 \times 10^{-3}$ & $1.61 \times 10^{-2}$ & $4.20 \times 10^{-8}$ & 384 & 147 016 320\\
9 & $1.95 \times 10^{-5}$ & 192 847    & 1.80 & $6.71 \times 10^{-4}$ & $8.14 \times 10^{-3}$ & $4.22 \times 10^{-8}$ & 768 & 148 106 496\\
10 & $9.75 \times 10^{-6}$ & 95 971    & 1.85 & $6.37 \times 10^{-4}$ & $4.01 \times 10^{-3}$ & $4.18 \times 10^{-8}$ & 1 536 & 147 411 456\\
11 & $4.88 \times 10^{-6}$ & 50 319    & 1.76 & $-1.48 \times 10^{-4}$& $2.18 \times 10^{-3}$ & $4.33 \times 10^{-8}$ & 3 072 & 154 579 968\\
12 & $2.44 \times 10^{-6}$ & 16 002    & 1.78 & $3.65 \times 10^{-4}$ & $4.87 \times 10^{-4}$ & $3.04 \times 10^{-8}$ & 6 144 & 98 316 288\\
13 & $1.22 \times 10^{-6}$ & 8 373     & 1.85 & $-1.06 \times 10^{-4}$& $1.75 \times 10^{-3}$ & $2.09 \times 10^{-7}$ & 12 288 & 102 887 424\\
14 & $6.10 \times 10^{-7}$ & 1 974     & 1.97 & $-6.83 \times 10^{-4}$& $1.78 \times 10^{-3}$ & $9.00 \times 10^{-7}$ & 24 576 & 48 513 024\\
15 & $3.05 \times 10^{-7}$ & 40        & 1.49 & $-1.05 \times 10^{-3}$& $1.49 \times 10^{-5}$ & $3.72 \times 10^{-7}$ & 49 152 & 1 966 080\\
\hline
$\sum$ & \multicolumn{4}{c}{} & & $1.94 \times 10^{-6}$ & & 1 742 633 519
\end{tabular}
}}
\end{table}

In these tables, we list the time step size $\Delta t_\ell$, number of samples $P_\ell$, variance of the fine simulations $\mathbb{V}\left[ F_\ell \right]$, expected value $\mathbb{E}\left[ F_\ell - F_{\ell-1} \right]$ and variance \change{$V_\ell$} of the differences of simulations, estimated variance of the estimator $\mathbb{V}[\hat{Y}_\ell]$, cost per sample $C_\ell$  and cost per level $P_\ell C_\ell$. The variance of the estimator at level $\ell$ \change{is estimated as} 
\begin{equation}
\mathbb{V}\left[ \hat{Y}_\ell \right ] = \frac{V_\ell}{P_\ell}.
\end{equation}
\add{
We see that the experimental results match the expected behavior of the multilevel Monte Carlo method. The number of samples $P_\ell$ needed to keep $\sum_{\ell=0}^L \mathbb{V}\left[ \hat{Y}_\ell \right ] < E^2$ decreases drastically in function of $\ell$. We also see that $\mathbb{E}\left[ F_L - F_{L-1} \right] < E^2$. The cost per level $P_\ell C_\ell$ is also spread quite evenly over the levels, once the time step is a couple orders of magnitude smaller than $\epsilon^2$. This is to be expected, as the geometric factor with which the cost increases with $\ell$ is asymptotically the same as that with which $V_\ell$ decreases. In short, we thus achieve the bias of the finest level, while a large amount of variance reduction is performed in the coarser levels.
}

The total cost of each multilevel simulation, relative to the cost of a single sample at the coarsest level is computed \change{as} the sum of the cost of each level. \change{We can estimate the cost for an equivalent classical Monte Carlo simulation by considering that one needs to perform}
\begin{equation}
P_C = \left\lceil \frac{\mathbb{V}\left[ F_L \right]}{\sum_{\ell=0}^L\mathbb{V}\left[ \hat{Y}_\ell \right ]} \right\rceil
\end{equation}
samples with the fine time step at level $L$, to achieve the same bias and variance as the multilevel estimator. The cost of each sample in the classic Monte Carlo estimator is $\frac{2}{3}C_L$, as we do not need to perform a correlated coarse simulation. Note that, for the numbers in \change{Table} \ref{mlmcaptab_summary}, $\mathbb{V}\left[ F_L \right]$ is estimated using very few samples, so these results should not be taken to literally. They do give the correct order of magnitude of the cost of the equivalent classical Monte Carlo method, however. We now compare the cost of the classical and multilevel Monte Carlo simulations in Table \ref{mlmcaptab_summary}.

\begin{table}
\caption{Cost comparison between classical and multilevel Monte Carlo\label{mlmcaptab_summary}}
\centering
\change{
\begin{tabular}{c | c | c | c}
\change{RMSE} & Classical cost & Multilevel cost & Speedup\\
\hline
0.1 & 4 544 & 8 062 & 0.56\\
0.01 & 28 627 968 & 10 102 066 & 2.83\\
0.001 & 25 167 200 256 & 1 742 633 519 & 14.4
\end{tabular}
}
\end{table}

As can be concluded from \change{Table} \ref{mlmcaptab_summary} the multilevel Monte Carlo method gives a significant computational advantage when we want to compute low bias results in the setting of the modified Goldstein-Taylor model. This speedup increases as the requested accuracy of the simulation is increased.

\add{
\subsubsection{Adding a coarse level}
\label{mlmcapsec_coarsexp}
It makes little sense to add a full sequence of levels in the regime where $\Delta t \gg \epsilon^2$. It could however make sense to add a single very coarse level to the simulation as the variance of $F_\ell$ is consistently larger than that of $F_\ell - F_{\ell-1}$ in Figures~\ref{mlmcapfig_exp1eps10} through~\ref{mlmcapfig_exp1eps001}. To test this idea we repeat the experiment as before, for $E = 0.01$, with a coarse level at $\Delta t_0=0.5$. The results of this experiment can be seen in Table~\ref{mlmcaptab_simulationresultscoarse}}

\begin{table}
\add{
\caption{Results of the simulation in Section \ref{mlmcapsec_expbenchmark} with an extra coarse level for $E=0.01$.\label{mlmcaptab_simulationresultscoarse}}
\centering 
\begin{tabular}{c | c | c | c | r | c | c | c || c}
Level & $\Delta t_\ell$ & $P_\ell$ & $\mathbb{V}\left[ F_\ell \right]$ & \multicolumn{1}{c|}{$\mathbb{E}\left[ F_\ell - F_{\ell-1} \right]$} & $V_\ell$ & $\mathbb{V}[\hat{Y}_\ell]$ & $C_\ell$ & $P_\ell C_\ell$ \\
\hline
0 & $5.00 \times 10^{-1}$ & 2 978 687 & 1.96 & $9.91 \times 10^{-1}$ & $1.96 \times 10^{0\phantom{-}}$ & $6.60 \times 10^{-7}$ & 0.02   & 59 574\\
1 & $1.00 \times 10^{-2}$ & 354 282   & 1.47 & $-1.26 \times 10^{-1}$ & $1.42 \times 10^{0\phantom{-}}$ & $4.01 \times 10^{-6}$ & 1.02   & 361 368\\
2 & $5.00 \times 10^{-3}$ & 114 863   & 1.47 & $1.06 \times 10^{-2}$ & $4.36 \times 10^{-1}$ & $3.79 \times 10^{-6}$ & 3   & 344 589\\
3 & $2.50 \times 10^{-3}$ & 77 905    & 1.57 & $3.14 \times 10^{-2}$ & $4.01 \times 10^{-1}$ & $5.15 \times 10^{-6}$ & 6   & 467 430\\
4 & $1.25 \times 10^{-3}$ & 47 439    & 1.68 & $2.92 \times 10^{-2}$ & $3.00 \times 10^{-1}$ & $6.32 \times 10^{-6}$ & 12  & 569 268\\
5 & $6.25 \times 10^{-4}$ & 27 466    & 1.78 & $2.57 \times 10^{-2}$ & $2.01 \times 10^{-1}$ & $7.32 \times 10^{-6}$ & 24  & 659 184\\
6 & $3.13 \times 10^{-4}$ & 14 599    & 1.75 & $9.28 \times 10^{-3}$ & $1.10 \times 10^{-1}$ & $7.56 \times 10^{-6}$ & 48  & 700 752\\
7 & $1.56 \times 10^{-4}$ & 7 666     & 1.71 & $3.99 \times 10^{-3}$ & $6.74 \times 10^{-2}$ & $8.79 \times 10^{-6}$ & 96  & 735 936\\
8 & $7.81 \times 10^{-5}$ & 3 195     & 2.04 & $4.04 \times 10^{-3}$ & $2.49 \times 10^{-2}$ & $7.80 \times 10^{-6}$ & 192 & 613 440\\
9 & $3.91 \times 10^{-5}$ & 40        & 1.54 & $5.51 \times 10^{-3}$ & $3.65 \times 10^{-3}$ & $9.12 \times 10^{-5}$  & 384 & 15 360\\
\hline
$\sum$ & \multicolumn{4}{c}{} & & $1.42 \times 10^{-4}$ & & 4 526 900
\end{tabular}
}
\end{table}

\add{
We see that the total cost of the simulation with the extra coarse level is lower than that of the simulation starting with $\Delta t=\epsilon^2$ (4 526 900 as apposed to 10 102 066). Based on this initial experiment, it makes sense to include a very coarse level when using the multilevel Monte Carlo method in this context. For a detailed analysis and more extensive numerical results, we refer to future work.
}

\section{Conclusion}
\label{mlmcapsec_conclusion}

In this work, we have derived a new multilevel scheme for asymptotic-preserving particle schemes of the form given in \eqref{mlmcapeq_GTmod}. We have demonstrated that this scheme has interesting convergence behavior as the time step is refined, which is apparent in the expected value and variance of sampled differences of the quantity of interest. On the one hand\add{,} we get the expected linear convergence to the exact model in terms of $\Delta t$ for a fixed value of $\epsilon$. On the other hand\add{,} we get convergence to pure diffusion in the limit for large values of $\Delta t$. This means that we can interpret the multilevel Monte Carlo method in this setting as refining upon an initial simulation of the heat equation, gradually including transport effects until the correct regime set by $\epsilon$ has been achieved. \change{We have shown that a significant speedup over classical Monte Carlo simulation is achieved when applying a geometric sequence of levels starting from $\Delta t = \epsilon^2$. We have also shown that adding an extra coarse level to the simulation further accelerates the computation in the considered test case.}

The approach taken in developing the asymptotic-preserving scheme is general\add{,} and it is straightforward to apply the coupling described in Section~\ref{mlmcapsec_correlation} to other, more general, models. As such, we are confident that the ideas expressed in this paper \change{will} also be applicable to more general equations than the Goldstein-Taylor model studied here. In future work, this scheme can, for example, be extended to higher dimensional models, both in terms of position and velocity. More complicated models including, for example, absorption terms can also be studied. \add{We intend to expand upon the results in Section~\ref{mlmcapsec_coarsexp}, as well as considering varying $\epsilon$ together with $\Delta t$.}

\begin{acknowledgement}
We thank Pieterjan Robbe for many helpful discussions on the multilevel Monte Carlo method. \add{We also thank the anonymous reviewers for their helpful suggestions for improving the quality of this work. The computational resources and services used in this work were provided by the VSC (Flemish Supercomputer Center), funded by the Research Foundation - Flanders (FWO) and the Flemish Government -- department EWI.}
\end{acknowledgement}

%
\bibliographystyle{spmpsci}
%

\end{document}